\documentclass{amsart}
\pdfoutput=1 
\usepackage{amsfonts,amssymb}
\usepackage{amsmath}
\usepackage[textwidth=15cm, textheight=22cm]{geometry}
\usepackage{amssymb}

\usepackage[ngerman, english]{babel}
\usepackage{babelbib}

\usepackage{bbm}

\usepackage{xspace}

\usepackage{booktabs}
\usepackage{adjustbox}

\usepackage{subfigure}

\usepackage{footnote}

\usepackage{multirow}	

\usepackage[table]{xcolor}

\usepackage[colorinlistoftodos,prependcaption,textsize=tiny]{todonotes}

\usepackage{epstopdf}

\usepackage{graphicx}%
\newtheorem{theorem}{Theorem}

\newtheorem{definition}[theorem]{Definition}

\newtheorem{proposition}[theorem]{Proposition}

\newcommand{\MM}{\mathcal{M}}

\newcommand{\E}{\mathbb{E}}

\newcommand{\N}{\mathbb{N}}

\newcommand{\R}{\mathbb{R}}

\newcommand{\ind}{\mathbbmss{1}}

\newcommand{\defgl}{\mathrel{\mathop:}=}
\newcommand{\defgr}{\mathrel={\mathop:}}

\DeclareMathOperator{\sign}{sgn}

\definecolor{verylightgray}{gray}{0.85}

\begin{document}
\title[The Euler scheme for SDEs with discontinuous drift coefficient]{The Euler scheme for stochastic differential equations with discontinuous drift coefficient: \\  A numerical study of the convergence rate}
\date{\today}
 \email[S. G\"ottlich]{} \email[S. G\"ottlich]{goettlich@uni-mannheim.de}
 \email[K. Lux]{} \email[K. Lux]{klux@mail.uni-mannheim.de}
\email[A. Neuenkirch]{neuenkirch@math.uni-mannheim.de} 
\author{S. G\"ottlich, K. Lux, A. Neuenkirch}
\address[S. G\"ottlich, K. Lux, A. Neuenkirch]{Department of Mathematics\\
University of Mannheim\\
68161 Mannheim, Germany}
\thanks{}

\begin{abstract}
	The Euler scheme is one of the standard schemes to obtain numerical approximations of stochastic differential equations (SDEs). 
	Its convergence properties are well-known in the case of globally Lipschitz continuous coefficients. However, in many situations, relevant systems do not show a smooth behavior, which results in SDE models with discontinuous drift coefficient.
	In this work, we will analyze the long time properties of the Euler scheme
	applied to SDEs with a piecewise constant drift and a constant diffusion coefficient and carry out intensive numerical tests for its convergence properties.
	We will emphasize on numerical convergence rates and analyze how they depend on properties of the drift coefficient and the initial value. We will also give theoretical interpretations of some of the arising phenomena.
	For application purposes, we will study a rank-based stock market model describing the evolution of the capital distribution within the market and provide theoretical as well as numerical results on the long time ranking behavior.
\end{abstract}

\maketitle




\section{Introduction}


In recent years, many applications related to stochastic differential equations (SDEs) with discontinuous drift coefficient
have emerged. These types of equations typically arise in mathematical finance and insurance \cite{Asmussen.1997,Ichiba.2013,Ichiba.2011,Karatzas.2009},
engineering applications~\cite{LaBolle.2000,Olama2009}, economy~\cite{Shardin.2016,Leobacher2015} or stochastic control problems \cite{Shardin.2016,Balakrish1980,kushner,touzi}.

The existence and uniqueness of solutions of SDEs in the standard case, i.e.\ the case of sufficiently smooth coefficients, is well understood~\cite{ks}. However, the standard theory on SDEs does not apply anymore in case of a discontinuous drift coefficient, e.g.\ a piecewise constant drift coefficient, and a special theory is 
needed to address the question of existence and uniqueness of solutions of such SDEs~\cite{Kleptsyna.1985,Veretennikov.1981,Zvonkin.1974}. 
The same is true for the numerical analysis: The convergence behavior of approximation schemes needs to be
reconsidered and ``research on numerical methods for SDEs with irregular coefficients is highly active.''\ (\cite[p.\ 2]{LSeuler}). In the case of a sufficiently smooth drift and a constant diffusion coefficient, the exact strong rate of convergence is 1 for the Euler scheme, see~\cite{DGR2003,Kloeden.1999}.
At the time, when the main part of the  research presented here was  undertaken,  no comparable result was known in the case of a discontinuous, e.g.\ piecewise constant, drift coefficient.
After many discussions and investigations, also inspired by a previous version of this manuscript, refined results are now about to be established, see Section \ref{subsec:equation}.

In this work, we will focus on numerical approximations of SDEs in the presence of a piecewise constant 
drift and a constant diffusion coefficient.  We will provide theoretical considerations on the long time behavior of approximated SDE solutions based on results from the theory of ergodic Markov chains. Moreover, we will provide further insight into the numerical behavior of  approximation schemes, in particular the Euler scheme, by analyzing the numerical convergence rates based on a reference solution.
The numerical speed of convergence heavily depends on
the initial value and properties of the drift coefficient, e.g. drift direction or
jump height. Our tests reveal that for a special class of drift coefficients the numerical convergence rates are higher and independent of
initial conditions due to the ergodicity of the Euler scheme and the underlying SDE. We also use the Euler scheme to verify the long time behavior of a rank-based stock market model~\cite{Banner.2005},
a prominent model in finance to describe the evolution of the capital distribution within the market. 

The remainder of this manuscript is as follows: 
In Section \ref{sec:Stability}, we will introduce some theoretical and numerical basics
and establish the ergodicity of the Euler approximations
in the case of an  appropriate, piecewise constant drift coefficient.
In Section \ref{sec:SimStudies}, we will discuss numerical convergence properties and further
findings of several numerical tests. 
We will conclude this work in Section \ref{sec:finance} with the  application from mathematical finance mentioned above, 
where SDEs of discontinuous type naturally arise. 


\section{Problem Description} \label{sec:Stability}


In this section, we will introduce our basic setting, i.e.\ the type of SDE, we are interested in, and some basic terms for the numerical tests. Besides the Euler scheme and its long time properties in our setting, we will also briefly discuss the applicability and performance of some other numerical schemes. 

\subsection{The Equation}  \label{subsec:equation}


In this manuscript, we will consider time-homogeneous SDEs with piecewise constant drift coefficient and additive noise:
\begin{align} \label{eq:SDE}
dX_t = \sum\limits_{j=1}^{s} \alpha_j \cdot \ind_{B_j}(X_t)dt + \sigma dW_t, \quad t \geq 0, \qquad X_0= \xi.
\end{align}
Here, we have $s \in \N$, $\alpha_j, \sigma, \xi \in \R$ and disjoint (possibly infinitely many) intervals $B_j \subset \mathbb{R}$   for all $1 \leq j \leq s$  and $\left( W_t \right)_{t \in [0,T]}$ is a one-dimensional Brownian motion.

The existence and uniqueness of solutions to this type of SDEs is guaranteed by results of \cite{Veretennikov.1981,Zvonkin.1974} and \cite{Kleptsyna.1985}. In \cite{Veretennikov.1981}, conditions on the drift and diffusion coefficient are derived under which the corresponding SDE has a unique strong solution. As emphasized therein, those conditions are in particular fulfilled for a bounded drift coefficient and a constant diffusion coefficient. Thus, the existence and uniqueness of a strong solution for SDEs of type \eqref{eq:SDE} is ensured.

\medskip

For the numerical analysis of SDEs with discontinuous drift and/or diffusion coefficient, the situation is more involved.
In this manuscript, we will focus on the strong convergence rate of the Euler scheme, which, for a general SDE
$$ dX_t = f(X_t) dt + g(X_t) dW_t, \qquad t \in [0,T],  \qquad X_0= \xi,$$ where $f$ and $g$ are such that a unique strong solution exists, is given by
\begin{align} \label{euler_alg}
	x_{k+1}^{\text{expE}} = x_k^{\text{expE}} + f(x_k^{\text{expE}}) \Delta + g(x_k^{\text{expE}}) (W_{(k+1)\Delta}-W_{k \Delta}), \quad k=0,\ldots,n-1, \qquad x_0^{\text{expE}}= \xi.
\end{align}
The underlying time discretization of the time interval $[0,T]$ is $0=t_0<t_1<\cdots<t_n=T$ with corresponding step size $\Delta\defgl \frac{T}{n}$, where $n+1$ is the number of grid points.

While its behavior is well-known for SDEs with Lipschitz continuous coefficients $f$ and $g$, much less has been known in more general cases, even for SDEs with additive noise and a piecewise constant drift coefficient. The first contribution in this area is -- up to the best of our knowledge -- the work
 \cite{Gyongy.1998}, where almost sure convergence of the Euler scheme has been established in the case of a one-sided Lipschitz drift coefficient,  a locally Lipschitz diffusion coefficient and the existence of a Lyapunov function for the SDE.
 The results of \cite{Halidias.2008}  give strong  convergence of the Euler scheme for  SDEs with additive noise in the case of a discontinuous, but monotone drift coefficient, while \cite{Simonsen2014} establishes
 the almost sure and strong convergence of the Euler scheme for SDEs with additive noise and drift of the form $f(x)=-\textrm{sign}(x)$.
 Recent contributions with respect to strong approximations of SDEs with discontinuous drift coefficient are a series of articles by Ngo and Taguchi \cite{NTbest,NTlast,NTmc} and  Leobacher and  Sz\"olgyeny \cite{LSeuler,LSbit,LSaap}, respectively. Very recently M\"uller-Gronbach and Yaroslavtseva \cite{MG-Y} established strong order $1/2$  for \eqref{euler_alg} in the case of a scalar equation with piecewise Lipschitz drift and non-additive noise and Neuenkirch et al.~obtained the same convergence order for an adaptive Euler scheme in the multi-dimensional case, see \cite{NSS}.
The weak approximation of SDEs with discontinuous coefficients has been studied in \cite{kohatsuhiga}, where an Euler-type scheme  based on an SDE with mollified drift coefficient is analyzed.

\medskip

In the case of SDE \eqref{eq:SDE}, the latest result on the strong convergence rate of the Euler scheme
 \begin{align} \label{eq:euler}
x_{k+1}^{\text{expE}} = x_k^{\text{expE}} + \sum\limits_{j=1}^{s} \alpha_j \cdot \ind_{B_j}(x_k^{\text{expE}}) \Delta + \sigma (W_{(k+1)\Delta}-W_{k \Delta}), \quad k = 0,\ldots,n-1, \quad x_0^{\text{expE}}= \xi,
\end{align}
for the approximation of $X_T$, i.e.\ the solution at time $T$, is an
    $L^2$-convergence order $3/4-\varepsilon$ in \cite{SN} for arbitrarily small $\varepsilon>0$.

For a better comparison, note that in the standard setting of an SDE with additive noise, where the drift coefficient is sufficiently smooth, the Euler scheme has an exact strong convergence order of $1$, see e.g. \cite{DGR2003} and \cite[p.\ 350f]{Kloeden.1999}.

So to summarize: The Euler scheme for our non-standard setting of SDE \eqref{eq:SDE} has at least $L^2$-convergence order $3/4-\varepsilon$. However, observing this convergence order numerically will be a different story (see Section \ref{sec:SimStudies}).

\subsection{Simulation studies and empirical convergence rates}
As already mentioned, we are interested in empirically measuring the strong convergence rate of the Euler scheme. The standard procedure for this is as follows:  The root mean squared  error (RMSE) at time $T$ for the Euler scheme \eqref{eq:SDE} with step size $\Delta =T/n$ is given by
\begin{align}
	e(n) \defgl  \left( \E\left \vert X_T - x_n^{\text{expE}} \right\vert^2 \right)^{1/2}. 
\end{align}
Since an explicit form  of $X_T$ is unknown in general, one needs to replace $X_T$ in our simulation studies by
a numerical reference solution $X_T^{\tt num}$, which is computed by the Euler scheme for an extremely small step size $\Delta=T/N$ with a very large number of $N+1$ grid points such that this approximation can be considered close enough to the true solution.
Moreover, also the expectation $ \E \vert X_T^{ \tt num} - x_n^{\text{expE}}\vert^2$ is not known explicitly, so we will approximate this expectation by the empirical RMSE
\begin{align}
e_{\tt emp}(n) = \sqrt{\frac{1}{M} \sum\limits_{i=1}^{M}\left\vert  \left( X_T^{\tt num} - x_n^{\text{expE}}\right)^{(i)}\right\vert^2},
\end{align}
with a large number $M$ of Monte Carlo repetitions, i.e.\ $( X_T^{\tt num} - x_n^{\text{expE}})^{(i)}$, $i=1, \ldots, M$, are iid copies of $X_T^{\tt num} - x_n^{\text{expE}}$.
Here, $X_T^{\tt num}$ and $x_n^{\text{expE}}$ have the same random input.
Note that $N$ has to be chosen sufficiently large to generate the numerical reference solution and to avoid oscillations  in $e_{\tt emp}(n)$, which might occur if $N$ and $n$ are close. The number of repetitions $M$ should also be large enough to have a good approximation of the expectation, i.e.\ a small Monte Carlo error.


\subsection{Other schemes}


A natural idea is of course to consider other schemes than the explicit Euler scheme and to compare them in our simulation studies.

\subsubsection{The implicit Euler scheme}

Implicit schemes have good stability properties, thus, they are a natural choice to consider. For an SDE with additive noise, where the drift coefficient is sufficiently smooth, the implicit Euler scheme has strong convergence order $1$ (see e.g. \cite{Alfonsi.2013} and \cite{Neuenkirch.2014}).

However, for SDEs of type \eqref{eq:SDE}, already the implicit Euler scheme is not well defined. To see this, consider the SDE
 \begin{align*} dX_t=  \left( \alpha_1  \cdot \ind_{(-\infty,0)} (X_t)  + \alpha_2  \cdot \ind_{[0,\infty)} (X_t)   \right) dt +  \sigma dW_t, \quad t \geq 0,  \qquad X_0=\xi, \end{align*} 
with $\alpha_1 >0 > \alpha_2$. The implicit Euler scheme 
$$ x_{k+1}^{\text{impE}} = x_{k}^{\text{impE}} +  \left(\alpha_1 \cdot \ind_{(-\infty,0)}(x_{k+1}^{\text{impE}}) +   \alpha_2 \cdot \ind_{[0,\infty)}(x_{k+1}^{\text{impE}})\right)\Delta + \sigma ( W_{(k+1)\Delta} - W_{k \Delta}), \ k=0,\ldots, n-1,$$
requires to solve, for fixed but arbitrary $z \in \mathbb{R}$,  the equation
$$ y - \left( \alpha_1  \cdot \ind_{(-\infty,0)} (y)  + \alpha_2  \cdot \ind_{[0,\infty)} (y)   \right) \Delta=  z,$$
with respect to $y \in  \mathbb{R}$. This equation  does not possess a solution if $ z \in (-\alpha_1 \Delta, -\alpha_2 \Delta)$, and hence an implicit Euler scheme is not well defined in this setting.

\subsubsection{The Heun scheme}

The Heun scheme is another scheme with strong order one for SDEs with additive noise under appropriate smoothness conditions on the drift coefficient. Adapted from \cite[p.\ 373]{Kloeden.1999} for SDEs of type \eqref{eq:SDE}, it is defined by 
\begin{align*} 
	x_{k+1}^{\text{Heun}} &= x_{k}^{\text{Heun}} + \frac{1}{2}  \left(\sum\limits_{j=1}^{s} \alpha_j \cdot \ind_{B_j}(x_{k}^{\text{Heun}}) + \sum\limits_{j=1}^{s} \alpha_j \cdot \ind_{B_j}(\Gamma_{k})\right) \Delta + \sigma (W_{(k+1)\Delta}-W_{k \Delta}), \\
	\Gamma_{k} &= x_{k}^{\text{Heun}} + \sum\limits_{j=1}^{s} \alpha_j \cdot \ind_{B_j}(x_{k}^{\text{Heun}}) \Delta + \sigma ( W_{(k+1)\Delta} - W_{k \Delta}), \qquad k=0,\ldots,n-1. \notag
\end{align*}

For a closer look at the behaviour of this scheme at a discontinuity assume that the drift coefficient is given by  $a(x)=\pm\textrm{sign}(x)$. An increment of the  Heun scheme  with $x_{k}^{\text{Heun}}=x$  is then given by
\begin{align*} 
	x_{k+1}^{\text{Heun}} - x =   \frac{1}{2} \big( a(x)  + a(x + a(x)\Delta  + \sigma   ( W_{(k+1)\Delta} - W_{k \Delta})) \big) \Delta 	
	+ \sigma ( W_{(k+1)\Delta} - W_{k \Delta}) .
\end{align*}
So if no drift change occurs in the Euler step $x + a(x)\Delta  + \sigma ( W_{(k+1)\Delta} - W_{k \Delta})$,  a Heun step and an Euler step coincide. However, if a drift change occurs in the Euler step, the Heun step
reads as
$$ x_{k+1}^{\text{Heun}} = x +    \sigma ( W_{(k+1)\Delta} - W_{k \Delta}), $$
i.e., it approximates the drift coefficient by zero and its dynamics are purely diffusion-based in this case. 

\subsubsection{A Wagner-Platen type scheme}

A strong order $1.5$-scheme for SDEs with smooth drift coeffi\-cient and additive noise is given by a Wagner-Platen type scheme (see e.g. \cite[p.\ 383]{Kloeden.1999}), which 
reads  in our setting as
\begin{align*}
	x_{k+1}^{\text{Pla}} = x_{k}^{\text{Pla}} & + a_{k} \Delta + \sigma  ( W_{(k+1)\Delta} - W_{k \Delta}) \\ &
	+ \frac{1}{4}  \left(a_{k}^{+} - 2a_k + a_{k}^{-}\right) \Delta 
	+ \frac{1}{2\sqrt{\Delta}} \left(a_{k}^{+} - a_{k}^{-}\right) \int_{k \Delta}^{(k+1)\Delta} (W_u-W_{k \Delta}) du,
\end{align*}
with
$$ 
\Gamma_{k}^{\pm} =x_{k}^{\text{Pla}} + a_k \Delta \pm \sigma \sqrt{\Delta}, \quad 
	a_k =a(x_{k}^{\text{Pla}}), \quad
	a_k^{\pm} = a(\Gamma_{k}^{\pm}), \quad \text{with} \ a(x) = \sum\limits_{j=1}^{s} \alpha_j \cdot \ind_{B_j}(x).$$
	
Now, we look again at the case of a drift coefficient  given by  $a(x)=\pm\textrm{sign}(x)$ and stepsize $\Delta < \sigma^2$. 
For a Wagner-Platen step with  $x_{k}^{\text{Pla}}=x$, it  depends  now on whether
$$  x + a(x) \Delta + \sigma \sqrt{\Delta}, \quad x, \quad x + a(x) \Delta - \sigma \sqrt{\Delta} $$
have the same sign or not. 
If this condition is fulfilled, i.e., if $x$ is sufficiently far away from the discontinuity, then, a Wagner-Platen step and an Euler step coincide.
If the latter condition is not satisfied, then we have the dynamics
\begin{align*}
	x_{k+1}^{\text{Pla}} = x &+  \frac{1}{2} a(x) \Delta + \sigma   ( W_{(k+1)\Delta} - W_{k \Delta}) \notag \\
	&+ \frac{1}{2\sqrt{\Delta}} \left(a\big(x+a(x) \Delta + \sigma \sqrt{\Delta}\big)  - a\big(x+a(x) \Delta - \sigma \sqrt{\Delta}\big)  \right) \int_{k \Delta}^{(k+1)\Delta} (W_u-W_{k \Delta}) du.
\end{align*}
So also here, the diffusive dynamic dominates the scheme when taking values close to the discontinuity.


\subsection{Ergodicity and stability of the Euler scheme} \label{sec:ergod}


We will now address the long time properties of the Euler scheme based on results from the theory of ergodic Markov chains. For simplicity, we consider here a special case of SDE \eqref{eq:SDE}, namely
\begin{align*} dX_t=  \left( \alpha_1  \cdot \ind_{(-\infty,0)} (X_t)  + \alpha_2  \cdot \ind_{[0,\infty)} (X_t)   \right) dt + dW_t, \quad t,s \geq 0,  \qquad X_0=\xi, \end{align*}
and assume that $$\alpha_1 >0 > \alpha_2,$$ i.e.\ a  drift coefficient, which is pointing towards zero.
Clearly, we have
\begin{align} \lim_{s \rightarrow 0} \mathbb{E}(X_{t+s}|X_t=x) = x +  \alpha_1  \cdot \ind_{(-\infty,0)} (x)  + \alpha_2  \cdot \ind_{[0,\infty)} (x), \qquad t \geq 0, \,\, x \neq 0, \label{ce_cont} \end{align}
i.e. on average, the solution is moving inwards. Moreover, following e.g.\ chapter 6 in \cite{gs}, this SDE admits a unique invariant distribution with Lebesgue density
$$ \varphi_{\infty}(x)= c \cdot e^{ 2 \alpha_2 x} \cdot \ind_{[0,\infty)} (x)  + c \cdot e^{ 2 \alpha_1 x}  \cdot \ind_{(-\infty,0)} (x), \qquad x \in \mathbb{R},$$
where the normalizing constant $c>0$ is such that $\int_{-\infty}^{\infty} \varphi_{\infty}(x) dx =1$.
In particular, we have
that \begin{align}\label{conv_cont}
 \lim_{t \rightarrow \infty} {\mathbb P}(X_t \leq y) = \int_{-\infty}^y \varphi_{\infty}(z) dz, \qquad y \in \mathbb{R}, \end{align}
and the law of large numbers
\begin{align} \lim_{L \rightarrow \infty} \frac{1}{L} \int_0^L h(X_t) dt =  \int_{-\infty}^{\infty} h(x) \varphi_{\infty}(x) dx  \qquad \textrm{a.s.,} \label{ln}	\end{align}
holds, if $h: \mathbb{R} \rightarrow \mathbb{R}$ is measurable and satisfies $ \int_{-\infty}^{\infty} |h(x)| \varphi_{\infty}(x) dx < \infty$.

\smallskip

It will turn out that the explicit Euler scheme
\begin{align}\label{e_disc} x_{k+1}^{\text{expE}, \xi}= x_k^{\text{expE}, \xi} + a(x_k^{\text{expE}, \xi}) \Delta + W_{(k+1)\Delta}-W_{k \Delta}, \quad k=0,1,\ldots,  \qquad x_0^{\text{expE}, \xi}= \xi,\end{align}
with
$$ a(x) = \alpha_1  \cdot \ind_{(-\infty,0)} (x)  + \alpha_2  \cdot \ind_{[0,\infty)} (x), \qquad x \in \mathbb{R},$$ will recover these properties. (Here we also indicate the dependence on the initial value $\xi$ in our notation.) The Euler scheme \eqref{e_disc}
corresponds to a time homogenous Markov chain with transition kernel
$$ p_{\Delta}(x,A) = \int_A \frac{1}{\sqrt{2 \pi \Delta}} \exp\left( -\frac{1}{2\Delta}\big(y-(x+a(x)\Delta \big)^2\right) dy, \qquad x \in \mathbb{R}, \quad A \in \mathcal{B}(\mathbb{R}),$$
and satisfies the discrete counterpart to \eqref{ce_cont}, i.e.\
\begin{align} \mathbb{E}( x_{k+1}^{\text{expE}, \xi}| x_{k}^{\text{expE}, \xi}=x) = x +  a(x)\Delta, \qquad k=0,1,\ldots, \qquad x \in \mathbb{R}. \label{ce_disc} \end{align}

Now, we will prove the existence of a unique stationary distribution for the Euler scheme. In particular, due to the discontinuity at zero, the following Proposition  \ref{ergod_prop} is not covered by the standard references as e.g.\ \cite{msh2002} and \cite{gr1996} for Euler-type discretizations of ergodic SDEs.
Note that the long time properties of \eqref{e_disc} have also been heuristically studied in \cite{Simonsen2013}.

However, we can easily verify that $V(x)=e^{ \tau |x|}$, $x \in \mathbb{R}$,  is an appropriate Lyapunov function for the above Markov chain, if $\tau>0$ is sufficiently small.
This is a direct consequence of the well known form of the moment generating function for the folded normal distribution, i.e.\ 
\begin{align}\label{folded} {\mathbb E} e^{\tau | \mu + \nu W_1|}=e^{\frac{\nu^2 \tau^2}{2}+\mu \tau}\left[1-\Phi\left(-\mu/\nu-\nu \tau \right) \right]+
e^{\frac{\nu^2 \tau^2}{2}-\mu \tau}\left[1-\Phi\left(\mu/\nu-\nu \tau \right) \right], \qquad \tau \in \mathbb{R}, \end{align}
where $\Phi$ is the distribution function of the standard normal distribution and $\mu \in \mathbb{R}$, $\nu >0$.
Using \eqref{folded}  with $\mu=x+a(x)\Delta$ and $\nu^2=\Delta$, we obtain
\begin{align*} {\mathbb{E}}\big( V (x_{k+1}^{\text{expE}, \xi}) | x_k^{\text{expE}, \xi} = x \big) & \leq e^{ \Delta \tau\left(\frac{\tau}{2}  + |\alpha_2|  \right)} +  e^{\Delta \tau \left( \frac{\tau}{2} - |\alpha_2| \right) }  e^{ \tau x}, \,\,\,\, \qquad x \geq 0, \\
{\mathbb{E}}\big( V (x_{k+1}^{\text{expE}, \xi}) | x_k^{\text{expE}, \xi} = x \big) & \leq e^{ \Delta \tau\left(\frac{\tau}{2}  + |\alpha_1| \right)} +  e^{\Delta \tau \left( \frac{\tau}{2} - |\alpha_1|\right) }  e^{-\tau x}, \qquad x < 0.
\end{align*}
So, we have
\begin{align*} {\mathbb{E}}\big( V (x_{k+1}^{\text{expE}, \xi}) | x_k^{\text{expE}, \xi} = x \big) & \leq e^{ \Delta \tau\left(\frac{\tau}{2}  + \max \{ |\alpha_1|,  |\alpha_2|  \}\right)} +  e^{\Delta \tau \left( \frac{\tau}{2} - \min\{ |\alpha_1|, |\alpha_2| \} \right) }  e^{ \tau |x|}, \qquad x \in \mathbb{R},
\end{align*}
and choosing $\tau < 2\min\{ |\alpha_1|, |\alpha_2| \}$ gives the desired property 
$$ {\mathbb{E}}\big( V (x_{k+1}^{\text{expE}, \xi}) | x_k^{\text{expE}, \xi} = x \big)  \leq C  +  \gamma  V(x), \qquad x \in \mathbb{R},
$$
with $C>0$, $\gamma \in (0,1)$. Since the transition kernel is Gaussian, an application of the quantitative Harris Theorem (see e.g. chapter 15 in \cite{Tweedie} or Theorem 3.15 (and the following example) in \cite{eberle}) yields the following geometric ergodicity result:

\begin{proposition}\label{ergod_prop}
Let $\alpha_1>0>\alpha_2$ and $\Delta>0$ be fixed. Then, the Euler scheme \eqref{e_disc} admits a unique stationary distribution $\mu_{\Delta}$, which is independent of the initial value $\xi$. Moreover, 
there exist $\beta_{\Delta} \in (0,1)$ and constants $\MM_{\Delta}(\xi), \xi \in \mathbb{R},$ such that
$$ \sup_{A \in \mathcal{B}(\mathbb{R})} \left| \mathbb{P}( x_{k}^{\textrm{expE}, \xi} \in A)  - \mu_{\Delta}(A) \right| \leq \MM_{\Delta}(\xi) \cdot \beta_{\Delta}^k, \qquad k \geq 1.$$
\end{proposition}

Choosing $A=(-\infty,y]$, we obtain in particular the counterpart to  \eqref{conv_cont}, i.e.\
\begin{align} \label{conv_disc}
 \lim_{k \rightarrow \infty} {\mathbb{P}}(x_k^{\text{expE}, \xi} \leq y) = \mu_{\Delta}((-\infty,y]), \qquad y \in \mathbb{R}. 
 \end{align}

Note that the limit distribution is independent of the initial value, as for the underlying SDE. 

\smallskip

Finally, an ergodic Theorem as e.g. Corollary 2.5 in \cite{eberle} yields also the discrete counterpart  to the law of large numbers \eqref{ln}: We have
\begin{align} \label{dln}
\lim_{L \rightarrow \infty}  \frac{1}{L}\sum_{k=1}^L h(x_k^{\text{expE}, \xi}) = \int_{-\infty}^{\infty} h(x) \mu_{\Delta}(dx)	\qquad \textrm{a.s.,}
\end{align}
 for all measurable $h: \mathbb{R} \rightarrow \mathbb{R}$ such that $\int_{-\infty}^{\infty} |h(x)|  \mu_{\Delta}(dx) < \infty.$

\bigskip


\section{Simulation Studies} \label{sec:SimStudies}


This section is concerned with the numerical investigation of SDEs of type \eqref{eq:SDE}. 
For the remainder, we will choose $T=1$, $M=10^5$, $N=2^{14}$ and $n=2^{\tilde{n}}$ with $\tilde{n} \in \{4,\ldots,10\}$ (unless otherwise mentioned).
We then calculate the corresponding Euler approximation and the empirical RMSE $e_{\tt emp}(n)$. For simplicity, we omit the upper index of the numerical approximation indicating that the approximation is based on the Euler scheme.
The empirical convergence rate is given by the negative slope of the regression line, which we obtain when plotting $\tilde{n}=\log_2(n)$ versus $\log_2\left(e_{\tt emp}(n)\right)$.
Here, we will focus on two types of drift coefficients: inward and outward pointing drift coefficients. 
\begin{definition}
We will call a drift coefficient $a: \mathbb{R} \rightarrow \mathbb{R}$ {\it inward pointing}, if there exists $x^* \in \mathbb{R}$ such that
$$ a(x)>0, \quad x < x^*, \qquad a(x)<0, \quad x>x^*,$$
and {\it outward pointing}, if there exists $x^* \in \mathbb{R}$ such that
$$ a(x)<0, \quad x < x^*, \qquad a(x)>0, \quad x>x^*.$$
\end{definition}

Our numerical investigations are based on several additional key characteristics:
We consider the average number of \textit{drift changes}. As the Euler scheme for SDE \eqref{eq:SDE} is exact up to the first drift change, another quantity of interest is the \textit{number of paths with at least one drift change}.
To get further insight whether some paths are really far away from the true solution, we measure the \textit{largest error} that occurs within the considered time interval (not necessarily in the end).
Besides the error sizes themselves, it is interesting to see what proportion of errors at final time $T$ is large, medium, or small and how this \textit{distribution of error sizes} depends on the step size.
Furthermore, we analyze the \textit{evolution of the error over time} for a fixed step size.
To underline the influence of the \textit{drift direction} towards or away from the discontinuity, we generate plots of several solution sample paths. We will see that the observed empirical\footnote{We use the expressions "numerical" and "empirical" rate (respectivley order) of convergence synonymously.} rates of convergence heavily depend on whether the drift coefficient is \textit{inward or outward pointing}. Whereas for the latter one, there is a dependency on the \textit{initial value} of the SDE, rates in case of an inward pointing drift coefficient seem to be independent of the initial value, corresponding to Proposition \ref{ergod_prop}.
In addition, we analyze how the \textit{jump height} (difference in drift values) influences the empirical convergence rate.

\medskip

As representatives of the class of SDEs \eqref{eq:SDE}, we consider here the SDEs given in Table \ref{tab:overviewSDEs}.

\begin{table}[htbp]
	\begin{minipage}{\textwidth}
		\centering
		\caption{Selection of analyzed SDEs}
				\label{tab:overviewSDEs}
		\begin{adjustbox}{max width=\textwidth}
			\begin{tabular}{ll}
				\textbf{Drift coefficient} & \textbf{Corresponding SDE} \\
				\hline sign & $dX_t = \sign(X_t)dt + dW_t$ \\ 
				\hline minusSign & $dX_t = -\sign(X_t)dt + dW_t$ \\
				\hline 10sign & $dX_t = 10 \cdot \sign(X_t)dt + dW_t$ \\ 
				\hline minus10sign & $dX_t = -10 \cdot \sign(X_t)dt + dW_t$ \\ 
				\hline elementary\_minus34 & $dX_t = \left(-3\cdot \ind_{(-\infty,1.4)}(X_t) + 4\cdot \ind_{[1.4,\infty)}(X_t)\right)dt + dW_t$ \\
				\hline elementary4minus3 & $dX_t = \left(4\cdot \ind_{(-\infty,1.4)}(X_t) - 3\cdot \ind_{[1.4,\infty)}(X_t)\right)dt + dW_t$ \\
				\hline elementary\_minus0.6\_1 & $dX_t = \left(-0.6\cdot \ind_{(-\infty,1.4)}(X_t) + \ind_{[1.4,\infty)}(X_t)\right)dt + dW_t$ \\
				\hline elementary1minus0.6 & $dX_t = \left(\ind_{(-\infty,1.4)}(X_t) - 0.6\cdot \ind_{[1.4,\infty)}(X_t)\right)dt + dW_t$ \\
				\hline 
			\end{tabular}
		\end{adjustbox}
	\end{minipage}
	\label{tab:OverviewAnalysedSDEs}
\end{table}

In the remainder of this chapter, we will present and discuss some key results of the simulation studies.


\subsection{Key results} \label{subsec:keyresults}


The  empirical convergence rates obtained by the Euler scheme are given in Table \ref{tab:EulerMSQRates-step4-finePart} (outward pointing drift coefficients highlighted in light gray, the discontinuity in gray): 
\begin{table}[htbp]
	\begin{minipage}{\textwidth}
	\setcounter{mpfootnote}{\value{footnote}}
	\renewcommand{\thempfootnote}{\arabic{mpfootnote}}
	\centering
	\caption{Numerical Euler convergence rates}
	\label{tab:EulerMSQRates-step4-finePart}%
	\begin{tabular}{rcccccc}
		\toprule
		\textbf{Initial values} & \textbf{-1} & \cellcolor{lightgray} \textbf{0} & \textbf{1} & \textbf{2.5} & \textbf{3} & \textbf{5} \\
		\midrule
		\rowcolor{verylightgray} sign  & 0.69 & \cellcolor{lightgray} 0.59 & 0.68 & 0.83 & 1.01 & \textbf{--} \\
		\rowcolor{verylightgray} 10sign & \textbf{--}\footnotemark[1] & \cellcolor{lightgray} 0.25  & \textbf{--} & \textbf{--} & \textbf{--} & \textbf{--} \\
		minusSign & 0.81 & \cellcolor{lightgray} 0.80 & 0.81 & 0.82 & 0.82 & 0.89 \\
		minus10sign & 0.91  & \cellcolor{lightgray} 0.91  & 0.91  & 0.91  & 0.91  & 0.91 \\
		\midrule
		\textbf{Initial values} & \textbf{0} & \textbf{1} & \textbf{1.2} & \textbf{1.25} & \cellcolor{lightgray} \textbf{1.4} & \textbf{2} \\
		\midrule
		\rowcolor{verylightgray} elementary\_minus34 & 1.17 & 0.37 & 0.38 & 0.40 & \cellcolor{lightgray} 0.39 & 0.31 \\
		\rowcolor{verylightgray} elementary\_minus0.6\_1 & 0.75  & 0.69  & 0.69  & 0.69  & \cellcolor{lightgray} 0.71  & 0.70 \\
		elementary4minus3 & 0.87 & 0.87 & 0.87 & 0.87 & \cellcolor{lightgray} 0.87 & 0.87 \\
		elementary1minus0.6 & 0.81  & 0.80  & 0.80  & 0.80  & \cellcolor{lightgray} 0.80  & 0.80 \\
		\bottomrule
	\end{tabular}%
	\footnotetext[1]{errors close to machine accuracy; no empirical convergence rate calculated (see also equation \eqref{exit_prob_2})}
	\setcounter{footnote}{\value{mpfootnote}}
	\end{minipage}
\end{table}%

Our results show that 
\begin{itemize}
 \item in general, we loose convergence order one, which the Euler scheme has under standard assumptions for SDEs with additive noise 
 \item and that a crucial factor is whether the drift coefficient is inward or outward pointing: for inward pointing coefficients the guaranteed convergence order $3/4$ is recovered, which is not always the case for outward pointing coefficients.
\end{itemize}
\medskip

Furthermore, our numerical tests show that neither using the Heun scheme  nor using the Platen scheme yields a different picture. In particular, convergence rates do not improve significantly, and the schemes do not yield a better resolution of the discontinuity (see Tables \ref{tab:HeunMSQRates-step4-finePart} and \ref{tab:PlatenMSQRates-step4-finePart}).

\begin{table}[htbp]
	\centering
	\caption{Numerical Heun convergence rates, step size $2^{-4}$ onwards}
	\begin{tabular}{rrrrrrr}
		\toprule
		\textbf{Initial values} & \textbf{0} & \textbf{1} & \textbf{1.2} & \textbf{1.25} & \textbf{1.4} & \textbf{2} \\
		\midrule
		elementary\_minus34 & 1.15 & 0.42 & 0.38 & 0.38 & 0.41 & 0.40 \\
		elementary4minus3 & 0.77 & 0.77 & 0.77 & 0.77  & 0.77 & 0.77 \\
		\bottomrule
	\end{tabular}%
	\label{tab:HeunMSQRates-step4-finePart}%
\end{table}%
\begin{table}[htbp]
	\centering
	\caption{Numerical Platen convergence rates, step size $2^{-4}$ onwards}
	\begin{tabular}{rcccccccc}
		\toprule
		\textbf{Initial values} & \textbf{0} & \textbf{1} & \textbf{1.2} & \textbf{1.25} & \textbf{1.4} & \textbf{2} \\
		\midrule
		elementary\_minus34 & 1.22 & 0.40 & 0.40 & 0.40 & 0.42 & 0.43 \\
		elementary4minus3 & 0.79 & 0.79 & 0.79 & 0.79 & 0.79 & 0.79 \\
		\bottomrule
	\end{tabular}%
	\label{tab:PlatenMSQRates-step4-finePart}%
\end{table}%


\subsection{Drift direction and initial value} \label{subsec:InwardOutward}


For an outward pointing drift coefficient, the numerical convergence order even seems to depend on the initial value and the spectrum of orders obtained for different initial values is very broad with values between $0.25$ and $1.17$ (see Table \ref{tab:EulerMSQRates-step4-finePart}). 

On the other hand, for an inward pointing drift coefficient, the convergence order seems to be independent of the initial value and the spectrum of orders numerically obtained for different initial values and inward pointing drift coefficients is tight with values between $0.80$ and $0.91$ (see Table \ref{tab:EulerMSQRates-step4-finePart}).
The stability of the estimates is due to the ergodicity of the SDE and the Euler scheme in this case, see Subsection \ref{sec:ergod}.
The geometric convergence speed in  Proposition \ref{ergod_prop}  explains 
why the numerical tests for inward pointing drift coefficients yield such stable estimates, independently of the initial value: $X_T^{\tt num}$ and $x_n$ are, for a sufficiently large number of grid points $n+1$, close to their unique
stationary distributions, which stabilizes the Monte-Carlo estimates. Also, as pointed out already above, the guaranteed convergence order $3/4$ is recovered here.

For the above equations, the structure of the drift coefficient is directly related to the number of drift changes. An inward pointing drift coefficient results in many drift changes, while in the case of an outward pointing drift coefficient, only few drift changes occur. We can further observe that:
\begin{itemize}
 \item[(i)] when starting away from the discontinuity, numerical rates for outward pointing drift coefficients are better than for inward ones;
 \item[(ii)] when starting close to the discontinuity, outward pointing drift coefficients imply worse numerical convergence rates than inward ones.
 \end{itemize}
So, in the latter case we obtain a positive correlation between the number of drift changes and the numerical convergence rate, which implies that frequent drift changes are not necessarily bad for the quality of the approximation -- quite the contrary seems to apply, which is surprising
at first glance.

\medskip

Hence, the type of monotonicity of the drift coefficient is of great importance.
Intuitively, an inward pointing drift coefficient should lead to many drift changes, which suggests that individual drift changes are not of great importance.
An outward pointing drift coefficient on the other hand pushes the solution away from the discontinuity implying a low number of drift changes.


\subsection{Jump height}


The intensity of the effects related to inward and outward pointing drift coefficients depends on the jump height, i.e.\ the distance between assigned drift values. In case of elementary\_minus34, this distance amounts to 7 whereas it is 1.6 in case of elementary\_minus0.6\_1. The empirical convergence rates in Table \ref{tab:EulerMSQRates-step4-finePart} show: The higher the jump height, the more pronounced are the effects described in Subsection \ref{subsec:InwardOutward}. Exemplary, there is a difference of $0.8$ in the empirical convergence rates for elementary\_minus34 for initial values $0$ and $1$ whereas this diffe\-rence is only $0.06$ for elementary\_minus0.6\_1.
This phenomenon is related to a scaling property. By enlarging the drift value, the influence of the diffusive part of the SDE is weakened:
Consider e.g. the SDE
$$ dX_t = \alpha \sign(X_t)dt + dW_t \label{eq:SDE X_t}$$
with $\alpha \geq 1$. Using the new variable $Y_t= \frac{1}{\alpha} X_t$ we have the dynamics
$$ dY_t=  \sign(Y_t)dt + \frac{1}{\alpha} dW_t,   $$
with a reduced diffusion  coefficient.



 
\subsection{Case study of an inward versus outward pointing drift coefficient}


In this Subsection, we will analyze the pattern described in \ref{subsec:InwardOutward} in more detail, exemplary for the drift coefficients elementary4minus3 and elementary\_minus34.

\subsubsection{Drift changes}
\begin{figure}[h!]
	\centering
	\includegraphics[width=0.6\linewidth]{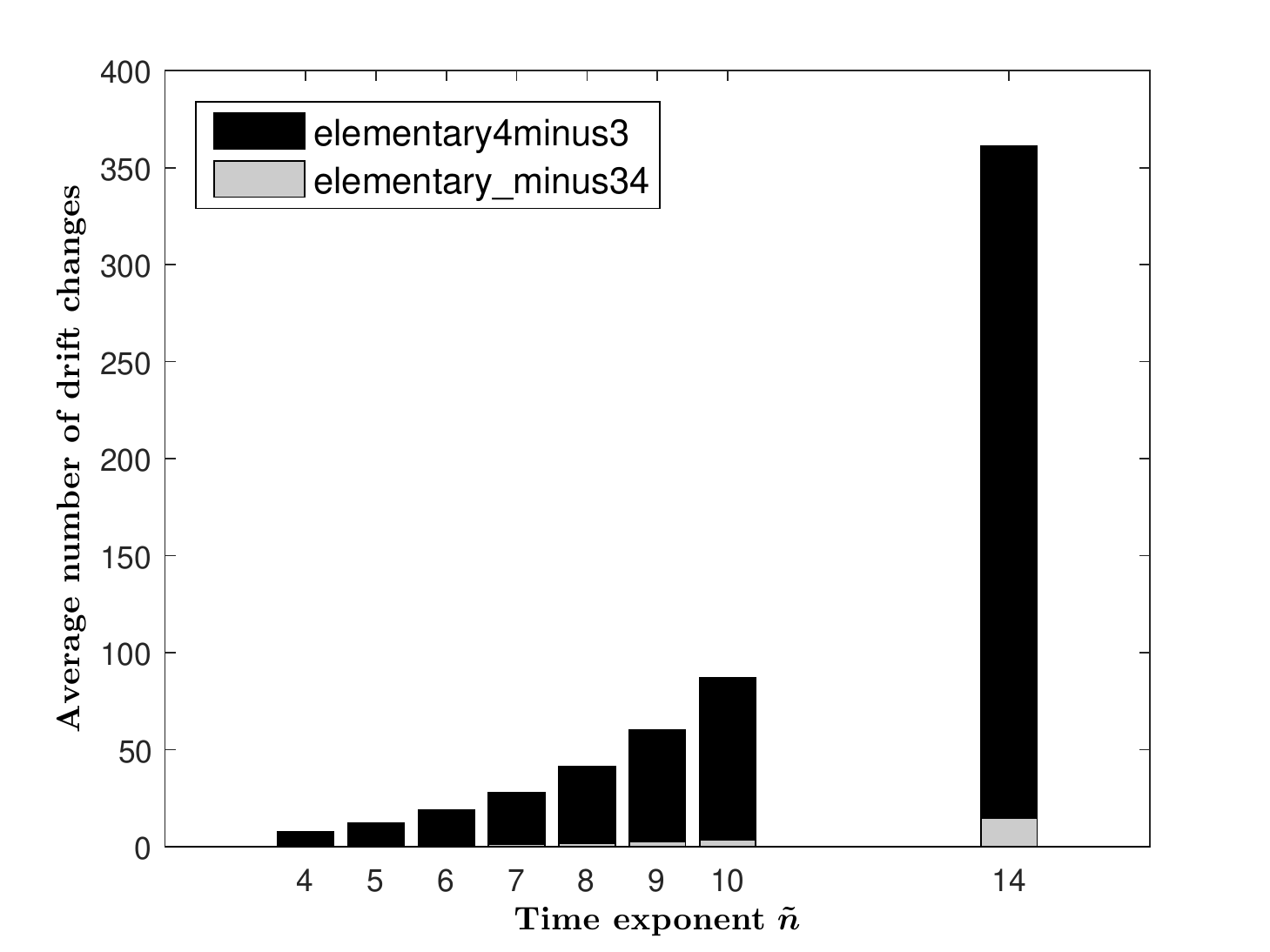}
	\caption{Average number of drift changes for $\xi=1.4$: elementary4minus3 vs. elementary\_minus34}
	\label{fig:avNumDriftChanges_ele4minus3ANDele_minus34_start1Komma4}
\end{figure}

Figure \ref{fig:avNumDriftChanges_ele4minus3ANDele_minus34_start1Komma4} shows the average number of drift changes for both coefficients.  The behavior goes along with the intuitive understanding described above.
Here  $\tilde{n}$ is the exponent of the dyadic step size $\Delta=2^{-\tilde{n}}$. Note that for step sizes $2^{-4}$ to $2^{-8}$ and elementary\_minus34 the number of drift changes stays below $2$.

\subsubsection{Comparison of solution sample paths}

\begin{figure}[h!]
	\subfigure[\ elementary4minus3, $\xi=1.4$]{\includegraphics[width=0.49\textwidth]{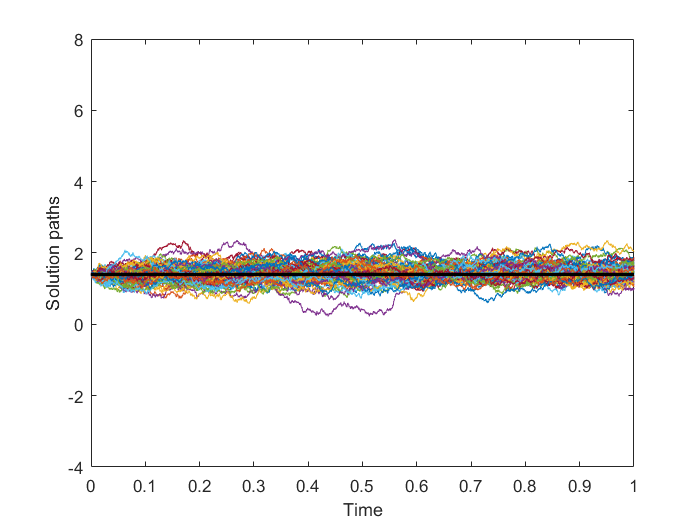}}\hfill
	\subfigure[\ elementary\_minus34, $\xi=1.4$]{\includegraphics[width=0.49\textwidth]{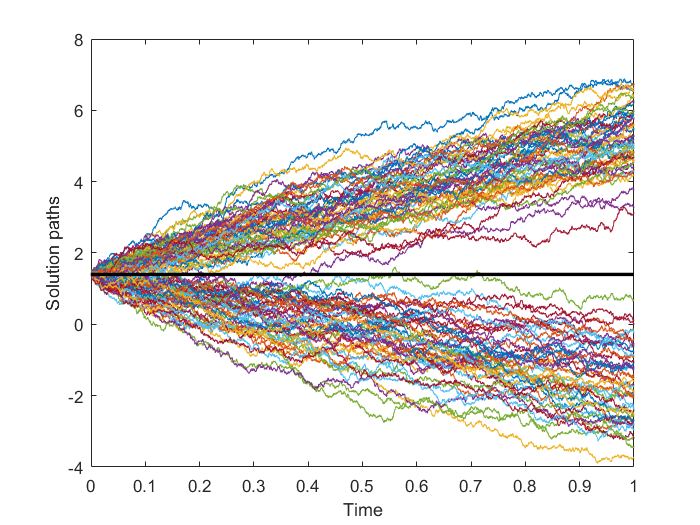}}

	\caption{Comparison of solution paths: elementary4minus3 vs. elementary\_minus34}
	
	\label{fig:samplePaths_ele4minus3VSele_minus34_start1Komma4}
\end{figure}

Figure  \ref{fig:samplePaths_ele4minus3VSele_minus34_start1Komma4} shows 100 sample paths of the numerical reference solution ($\Delta=2^{-14}$). The black line represents the discontinuity in the drift coefficient.

\smallskip

\begin{figure}[h!]
	\centering
	\includegraphics[width=0.6\linewidth]{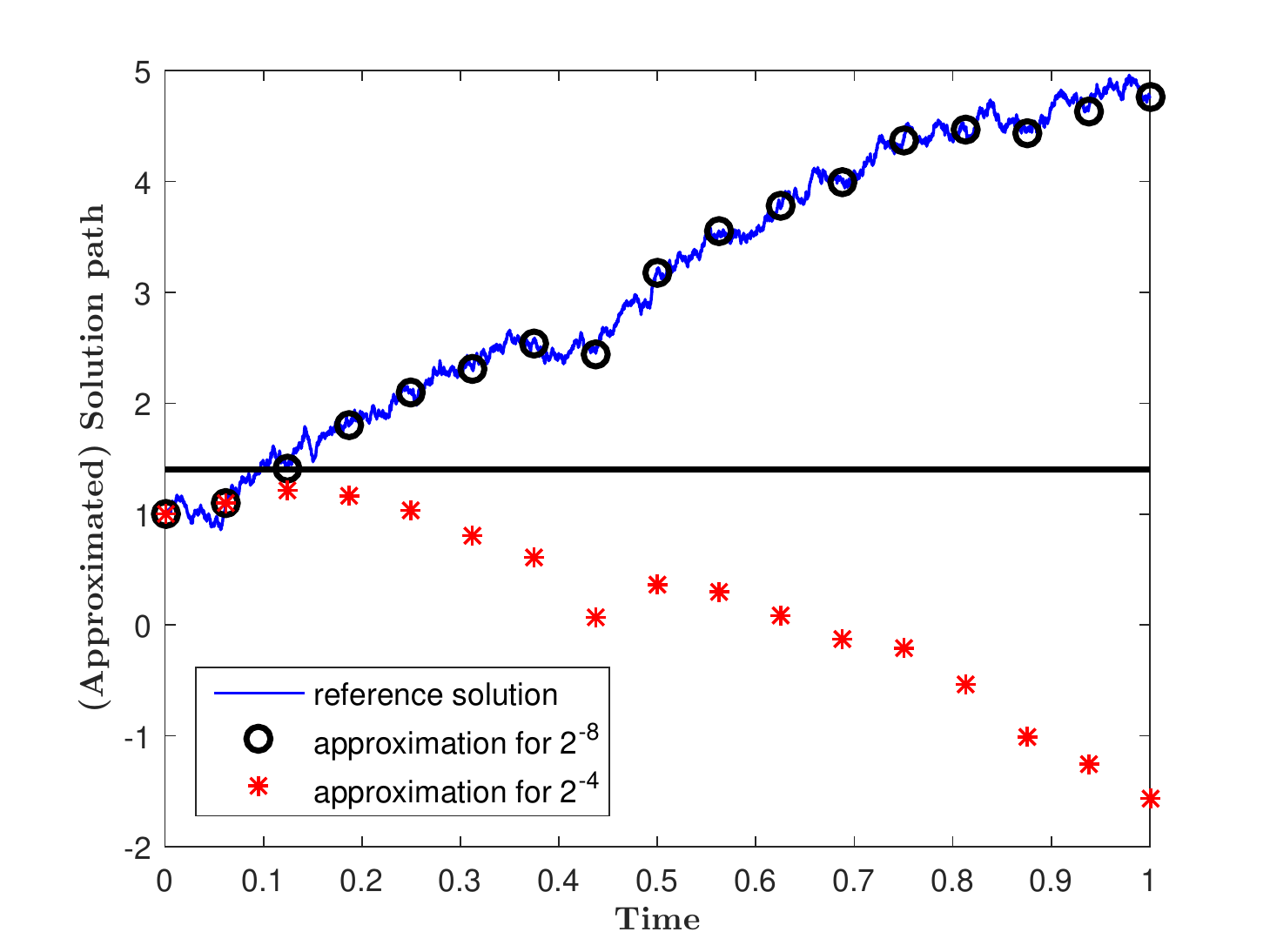}
	\caption{Importance of capturing the drift changes for elementary\_minus34, $\xi=1$}
	\label{fig:importance_captureDriftChange_ele_minus34_1_sample9999}
\end{figure}

In the situation of Figure \ref{fig:samplePaths_ele4minus3VSele_minus34_start1Komma4}(b), where the solution drifts away from the discontinuity, it is of tremendous importance whether a drift change is captured by the approximation or not: the solution does not stay close to the discontinuity and thus, there are not many chances for a drift correction to take place, see Figure \ref{fig:importance_captureDriftChange_ele_minus34_1_sample9999}.
For the SDE
\begin{align} dX_t=  \left( \alpha_1  \cdot \ind_{(-\infty,0)} (X_t)  + \alpha_2  \cdot \ind_{[0,\infty)} (X_t)   \right) dt + dW_t, \quad t \geq 0,  \qquad X_0=\xi, \label{test_illustrate} \end{align}
with $\alpha_1<0<\alpha_2$ and $\xi >0$ the conditional probability $p(\xi,\theta,\Delta)$ that the exact solution changes its drift over $[0,\Delta]$ given that the approximation $x_1$ at $t=\Delta$ has value $\theta \geq 0$ (and thus has not changed its drift) satisfies
\begin{align} p(\xi,\theta,\Delta):=\mathbb{P}\Big{(} \inf_{t \in [0,\Delta]}X_t<0 \Big{|} X_0=\xi, x_1 = \theta \Big{)}= \exp \left( - 2 \frac{\xi \theta}{\Delta} \right), \end{align}
see e.g. \cite{Gobet}, page 169. So the (conditional) probability of missing drift changes is not negligible and even close to one for small $\xi$  or $\theta$.

\subsubsection{Largest error}

The latter observation is also reflected in the largest distance for $10^4$ sample paths between the approximation based on step size $2^{-10}$ and the numerical reference solution, see Table \ref{tab:errorSizes_MSQ_Euler}.
The largest distances amount to 1.271 for elementary4minus3 and 4.508 for elementary\_minus34. 
\begin{table}[htbp]
	\centering
	\caption{Largest and smallest Euler errors}
	\begin{tabular}{rrrrrrrrrr}
		\toprule
		\textbf{Initial values} & \textbf{} & \textbf{0} & \textbf{1} & \textbf{1.2} & \textbf{1.25} & \textbf{1.4} & \textbf{2} \\
		\midrule
		\textbf{elementary\_minus34} & max   & 0.045 & 1.331 & 2.614 & 3.087 & 4.508 & 0.335 \\
		& min   & 0.0002 & 0.179 & 0.331 & 0.383 & 0.696 & 0.059 \\
		\textbf{elementary4minus3} & max   & 1.223 & 0.934 & 0.968 & 0.981 & 1.013 & 1.271 \\
		& min   & 0.005 & 0.005 & 0.005 & 0.005 & 0.005 & 0.005 \\
		\bottomrule
	\end{tabular}%
	\label{tab:errorSizes_MSQ_Euler}%
\end{table}%


\subsubsection{Evolution of the error over time}

\begin{figure}[htbp]
	\subfigure[\ elementary4minus3, $\tilde{n}=4$]{\includegraphics[width=0.45\textwidth]{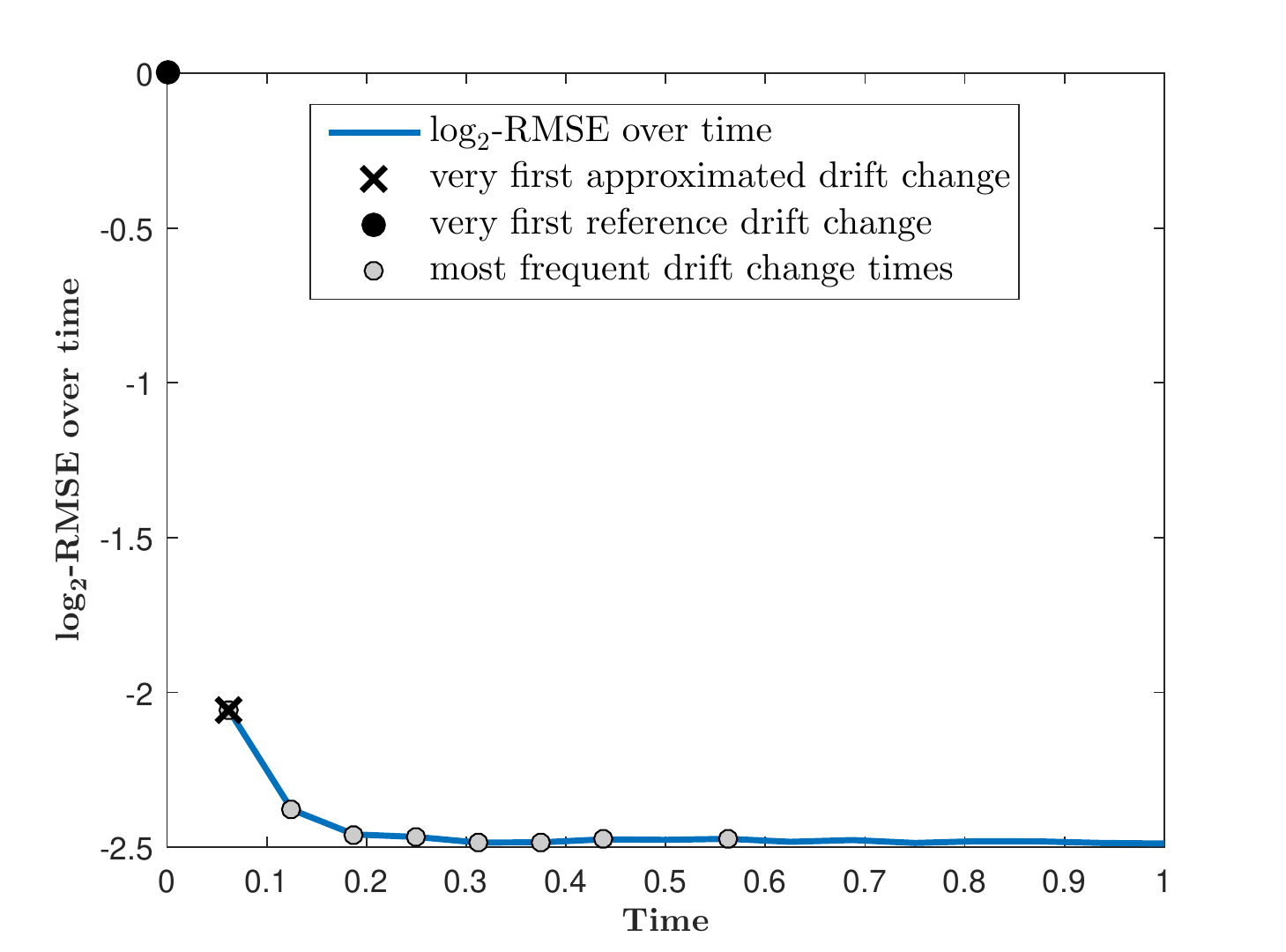}}
	\subfigure[\ elementary\_minus34, $\tilde{n}=4$]{\includegraphics[width=0.45\textwidth]{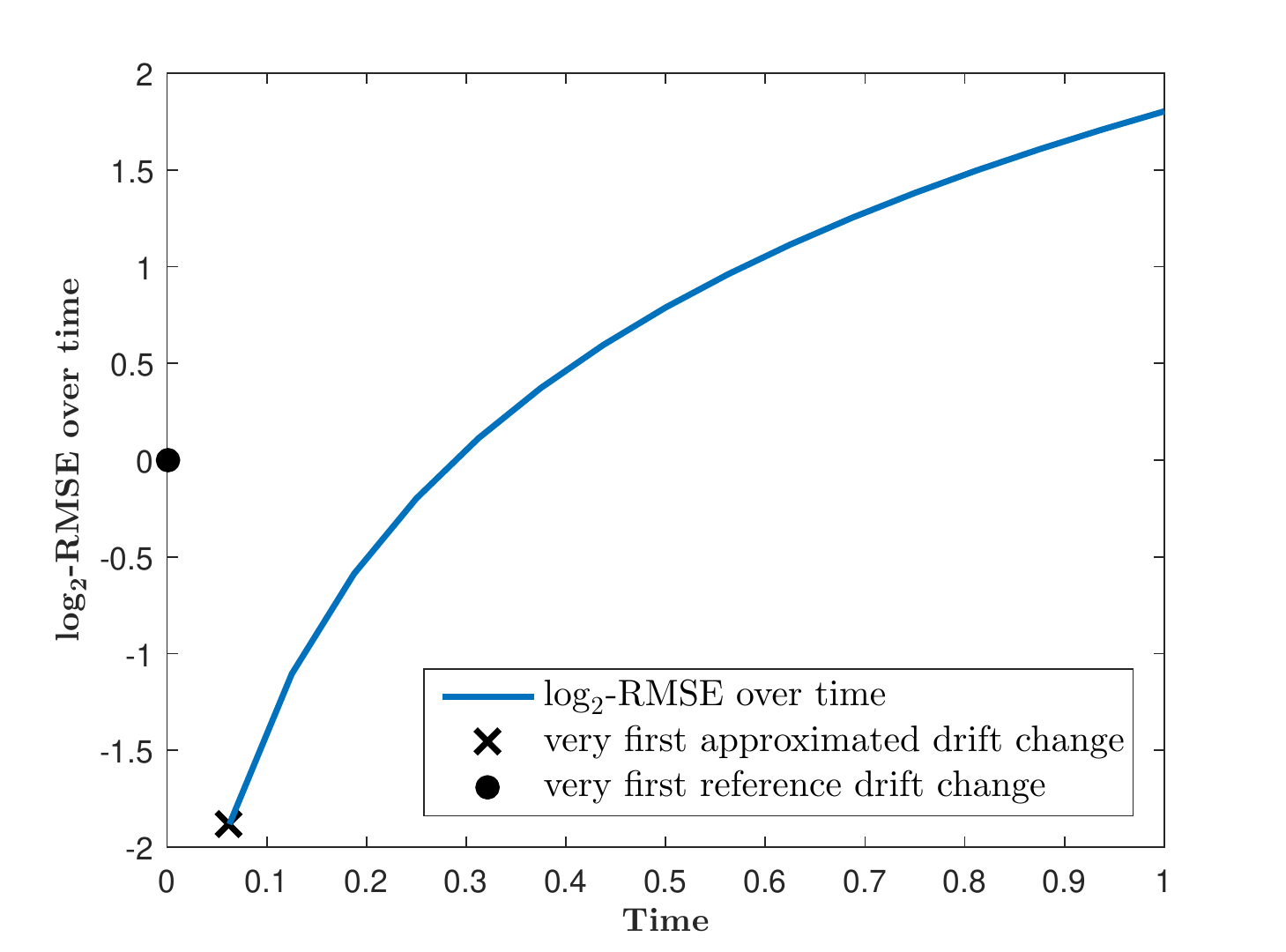}}
	\subfigure[\ elementary4minus3, $\tilde{n}=8$]{\includegraphics[width=0.45\textwidth]{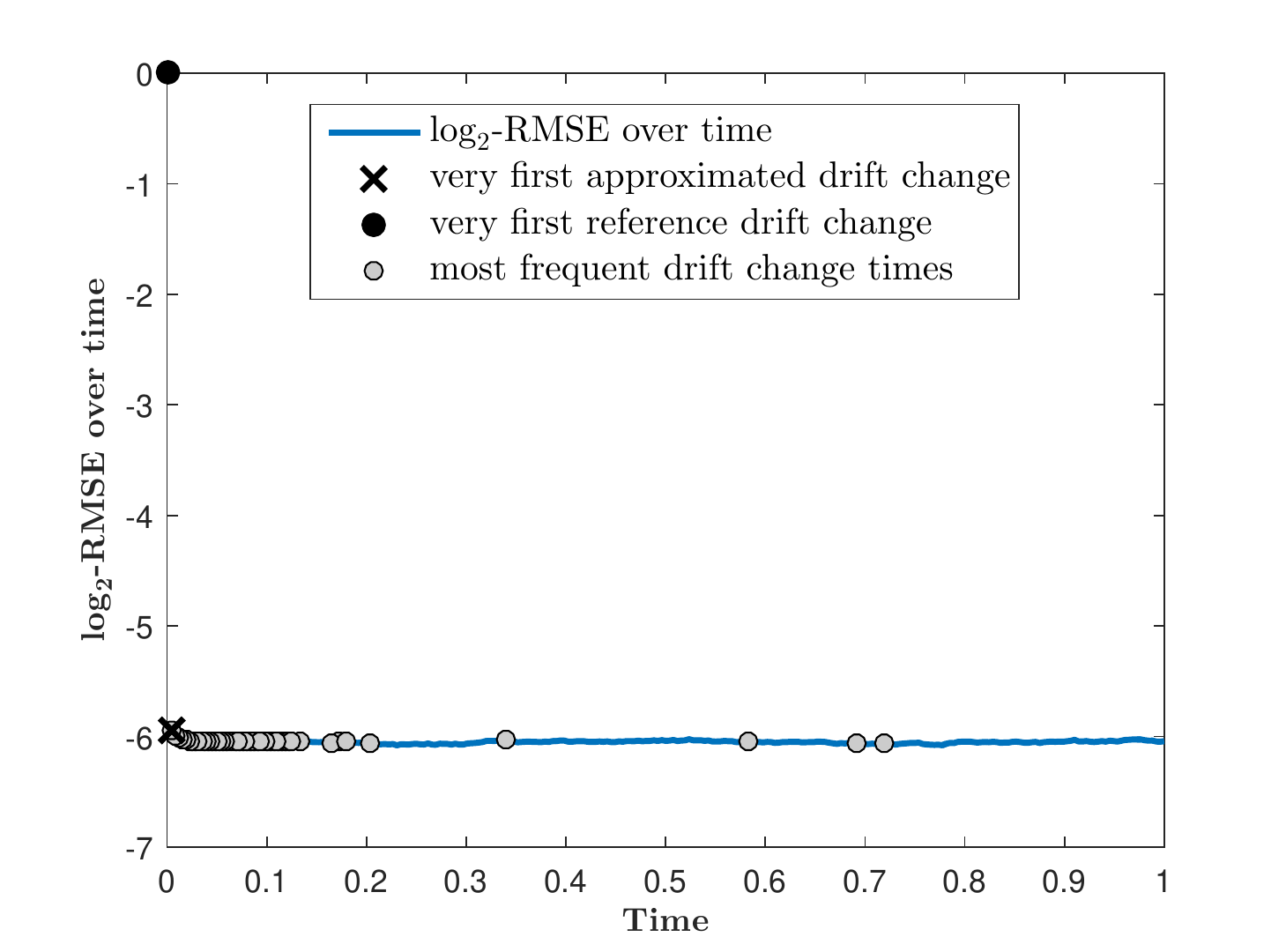}}
	\subfigure[\ elementary\_minus34, $\tilde{n}=8$]{\includegraphics[width=0.45\textwidth]{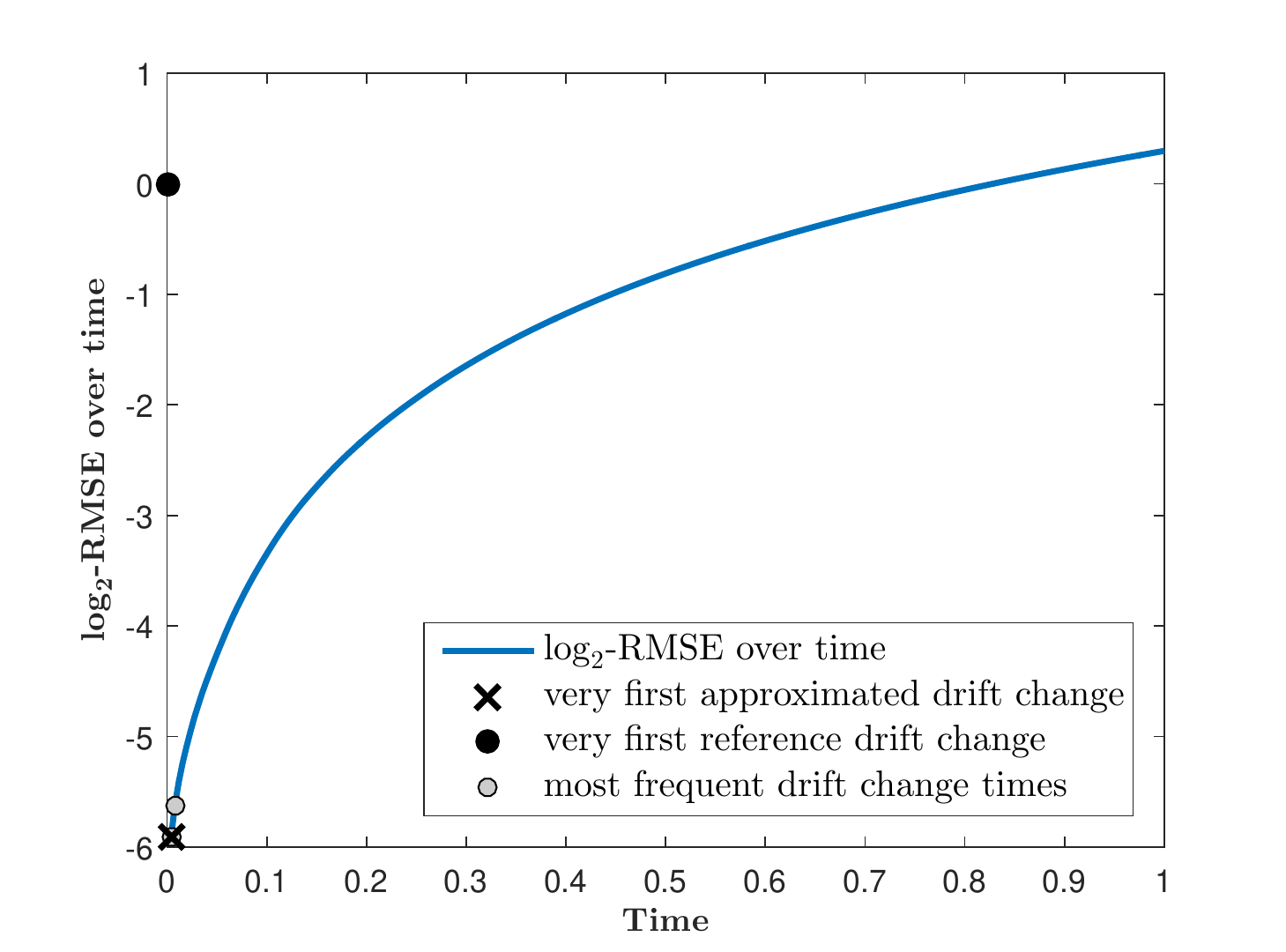}}
	\subfigure[\ elementary4minus3, $\tilde{n}=10$]{\includegraphics[width=0.45\textwidth]{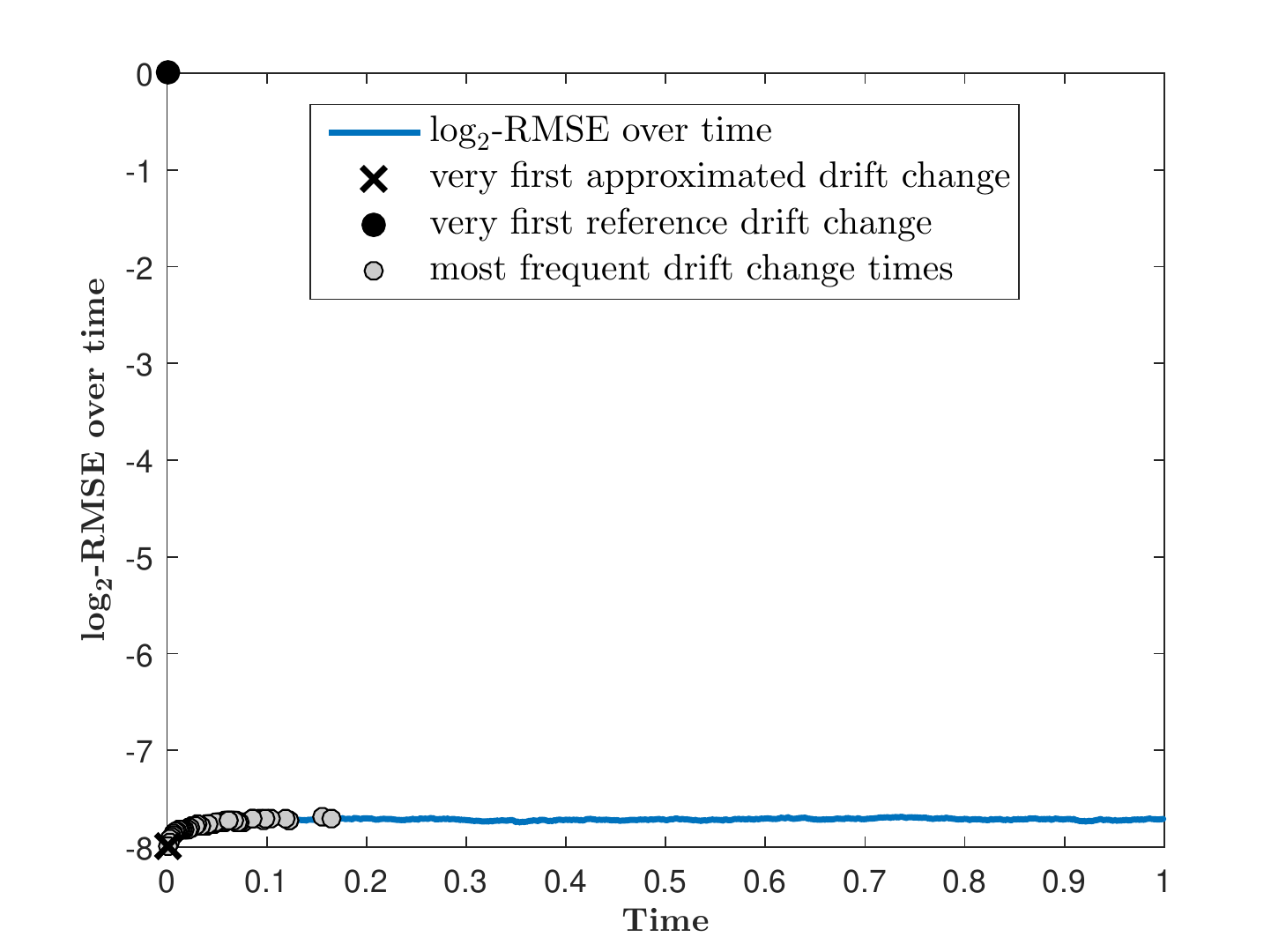}}
	\subfigure[\ elementary\_minus34, $\tilde{n}=10$]{\includegraphics[width=0.45\textwidth]{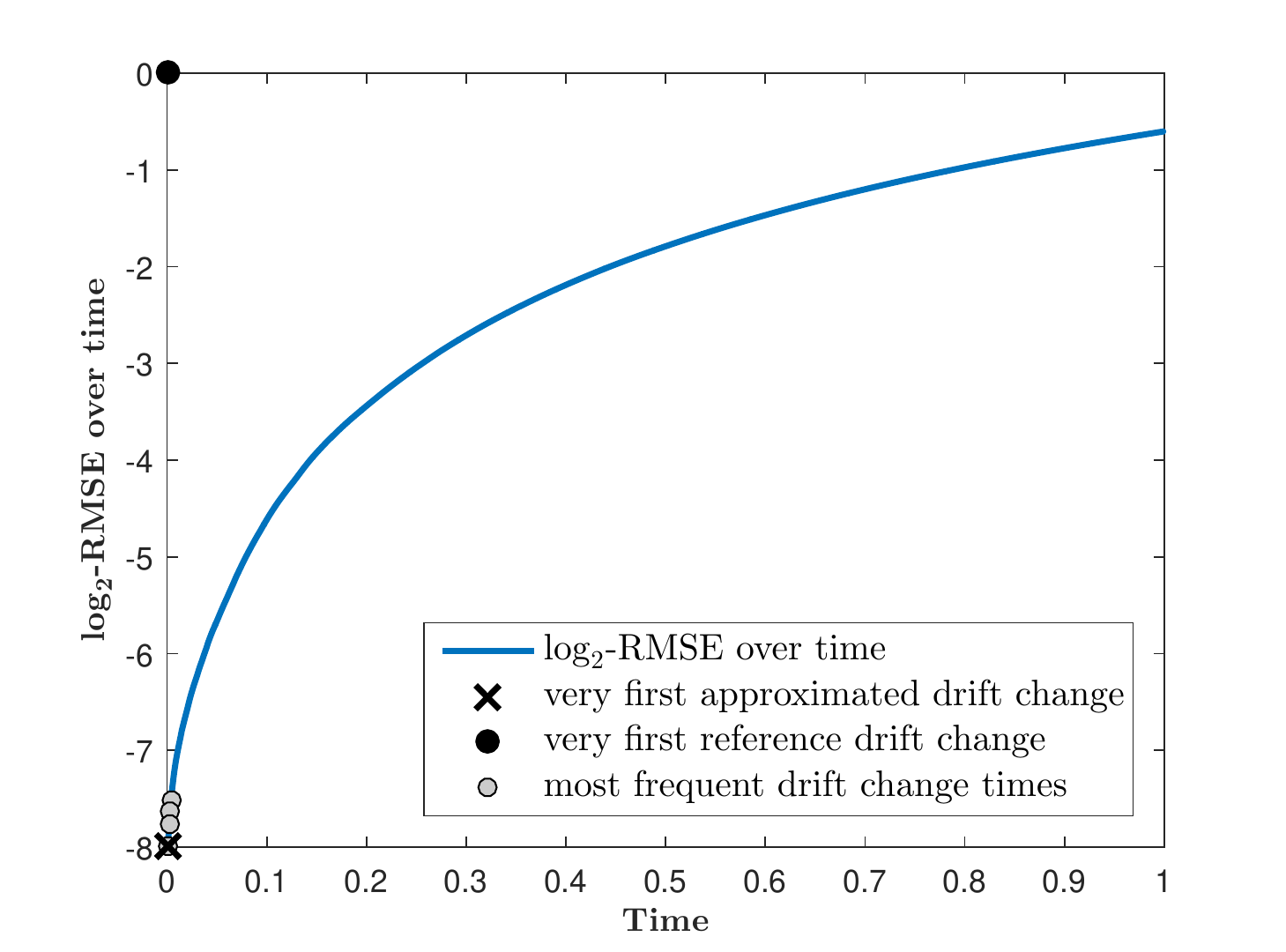}}
	
	\caption{Comparison of the error evolution over time for $\xi=1.4$ for different step sizes: elementary4minus3 vs. elementary\_minus34}
	
	\label{fig:errorOverTime_4_8_10_ele4minus3ANDele_minus34_start1Komma4}
\end{figure}

To gain even more insight, we compare the empirical RMSE for increasing time $t$ 
of elementary4minus3 and elementary\_minus34 when starting in the discontinuity $\xi=1.4$ for step sizes $2^{-4}$, $2^{-8}$ and $2^{-10}$ by plotting the base-2 logarithm of the RMSE against the time (see Figure \ref{fig:errorOverTime_4_8_10_ele4minus3ANDele_minus34_start1Komma4}).
We have added in these figures the following additional information: If the number is not zero, the most frequent times of drift changes corresponding to the chosen step size are indicated. The number of plotted drift change times is based on the average number of drift changes over the simulated sample paths.

Furthermore, if in the corresponding cases drift changes occur, we add the very first drift change (of all simulated paths) of the numerical reference solution and the Euler schemes. They are generated by finding the time at which the first drift change occurs for $10^4$ saved paths and then taking the minimum over all that times. The time is registered as the point of discretization at which a drift change that took place was detected. The very first drift change of the reference solution is marked at a height of zero for a better distinguishability. 
RMSE over time and drift change times are calculated on a basis of $10^{4}$ simulation paths.

We can extract from Figure \ref{fig:errorOverTime_4_8_10_ele4minus3ANDele_minus34_start1Komma4} at least two features:
\begin{itemize}
 \item[(i)] The error stays constant or even decreases over time for elementary4minus3 -- in contrast to a strong error accumulation over time for elementary\_minus34. (Note that the ordinate has a base-2 log scale.)
 \item[(ii)] In the inward pointing drift coefficient case, the error is by several magnitudes smaller than for an outward pointing drift coefficient. 
\end{itemize}

This illustrates again the  stabilizing effect of an inward pointing drift coefficient and the importance of capturing the first drift changes correctly in case of an outward pointing drift coefficient.

\subsubsection{Distribution of error sizes}

\begin{figure}[htbp]
	\subfigure[\ elementary4minus3]{\includegraphics[width=0.49\textwidth]{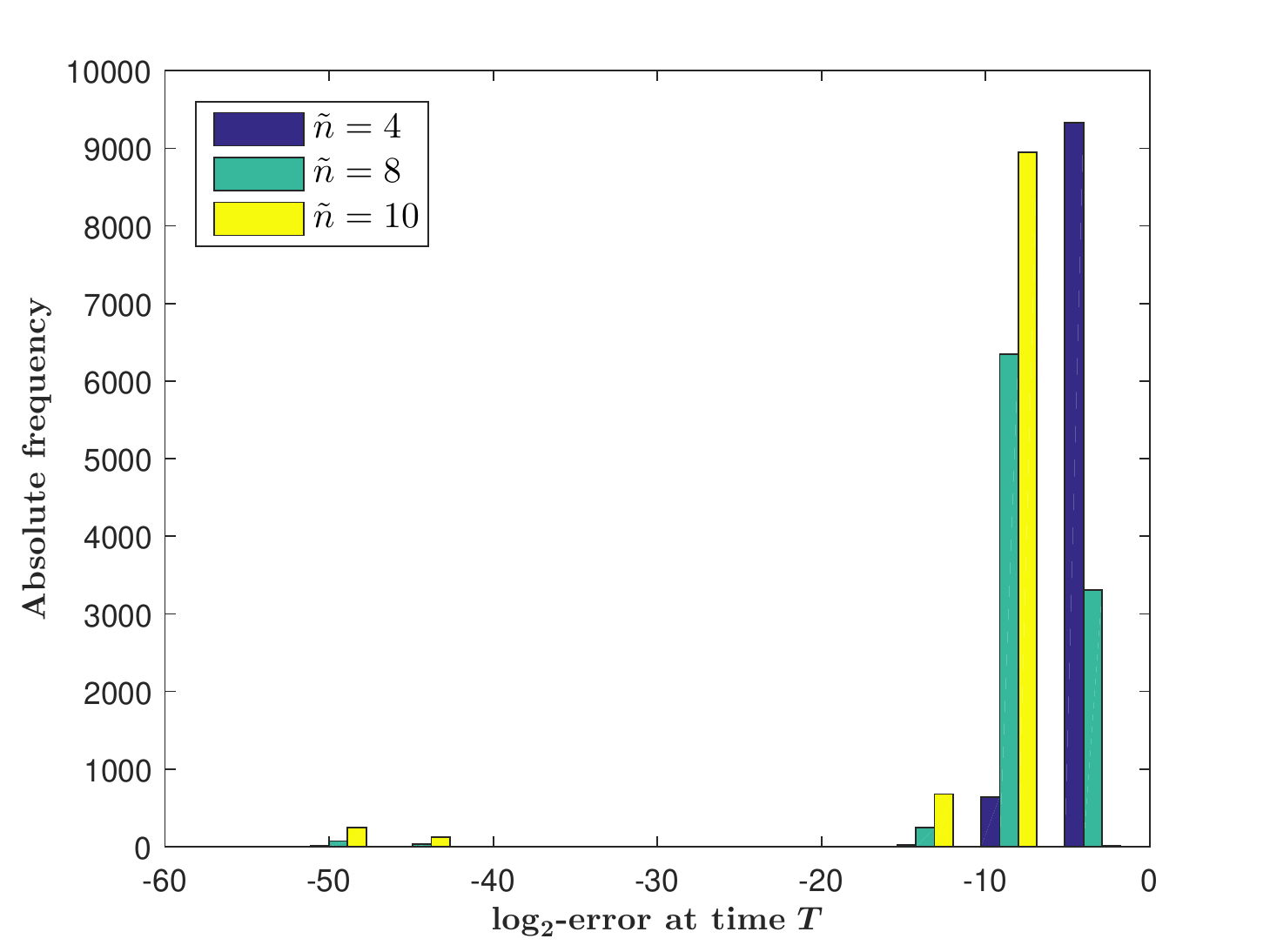}}\hfill
	\subfigure[\ elementary\_minus34]{\includegraphics[width=0.49\textwidth]{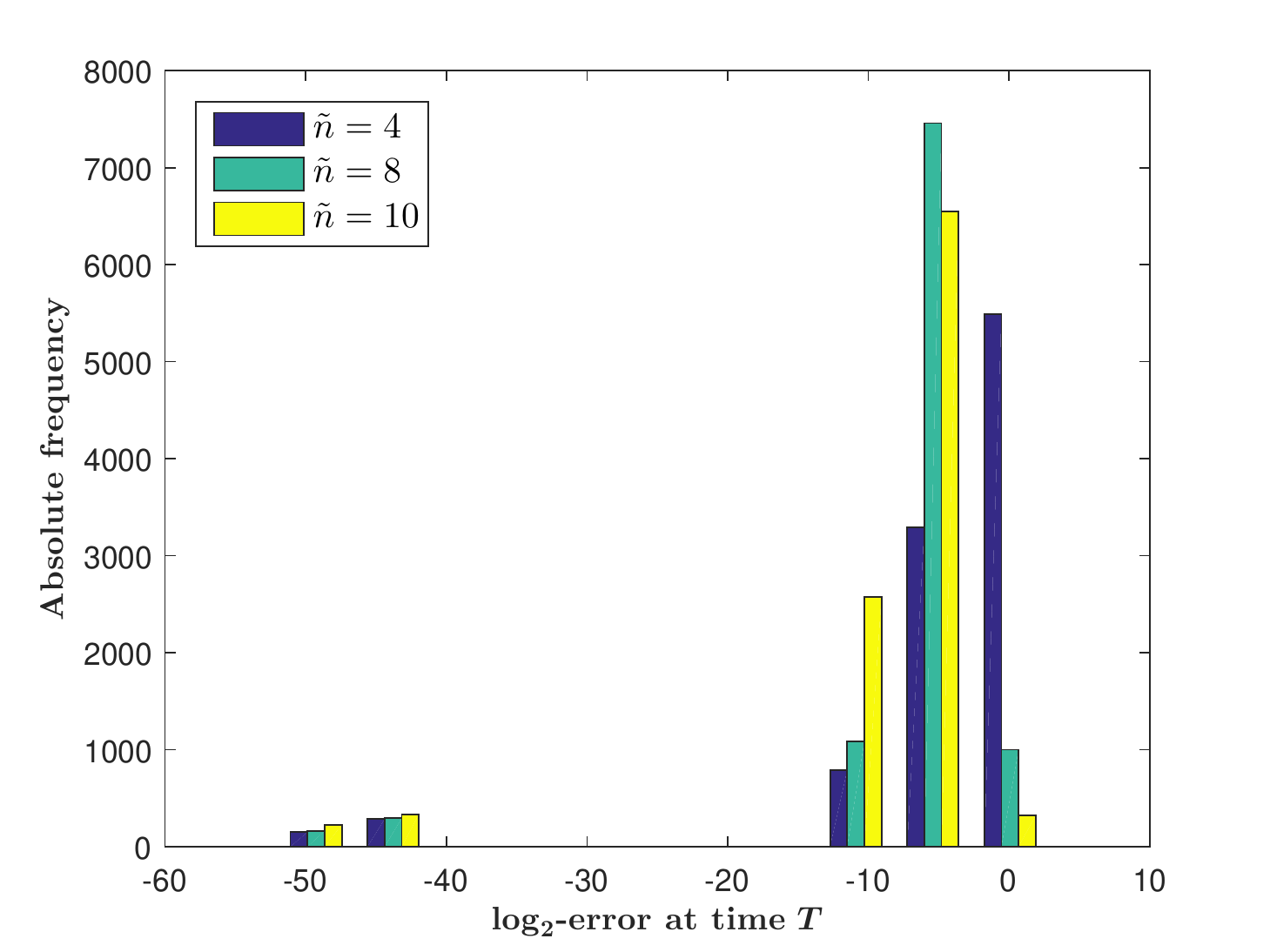}}
	
	\caption{Distribution of the error at time $T$ for different step sizes for elementary4minus3 and elementary\_minus34 with $\xi=1.4$}
	
	\label{fig:hist_10power5_ele4minus3_ele_minus34_n4810_1Komma4}
\end{figure}

Besides the empirical RMSE itself, the empirical distribution of the errors in $t=T$ is of interest. The error at final time $T$ is quantified by $\left\vert x_N - x_n\right\vert$ for step size $\Delta=T/n=2^{-\tilde{n}}$. The histograms in Figure \ref{fig:hist_10power5_ele4minus3_ele_minus34_n4810_1Komma4}  are based on $M=10^4$ simulations for different step sizes and highlight again the different magnitudes of the  empirical RMSE (abscissa with a base-$2$ logarithm scale).
Another feature, which we can extract from the histograms, is a non-negligible part of simulated paths with an error of machine accuracy size for elementary\_minus34. We will discuss this feature in more detail in the next Subsection. 

\subsection{Rare events and goodness of the regression fit}

In case of an outward pointing drift coefficient the empirical RMSE and the linear regression estimates become unreliable or at least questionable.

\begin{figure}[htbp]
	\subfigure[\ $\xi=0$]{\includegraphics[width=0.49\textwidth]{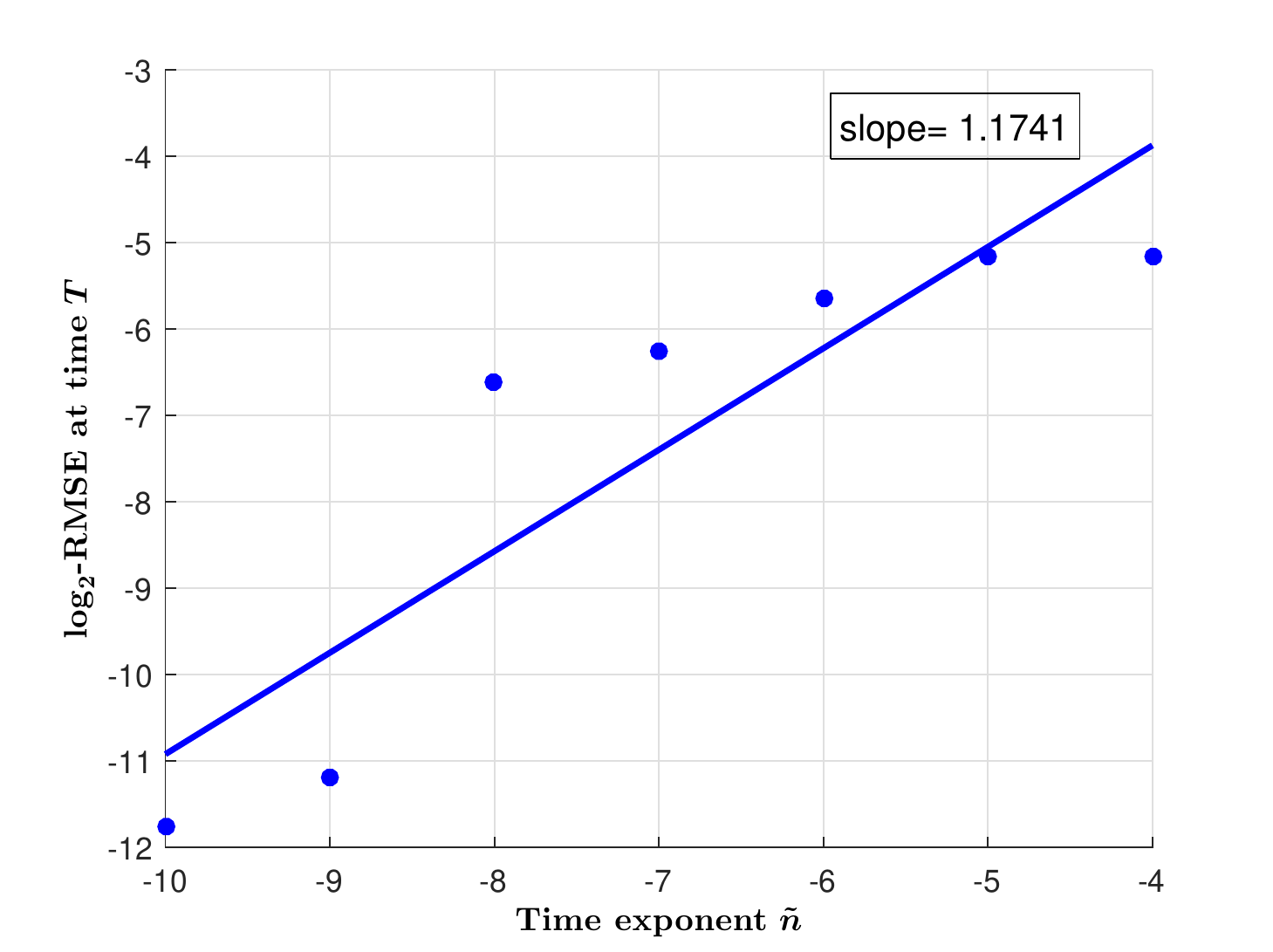}}\hfill
	\subfigure[\ $\xi=1$]{\includegraphics[width=0.49\textwidth]{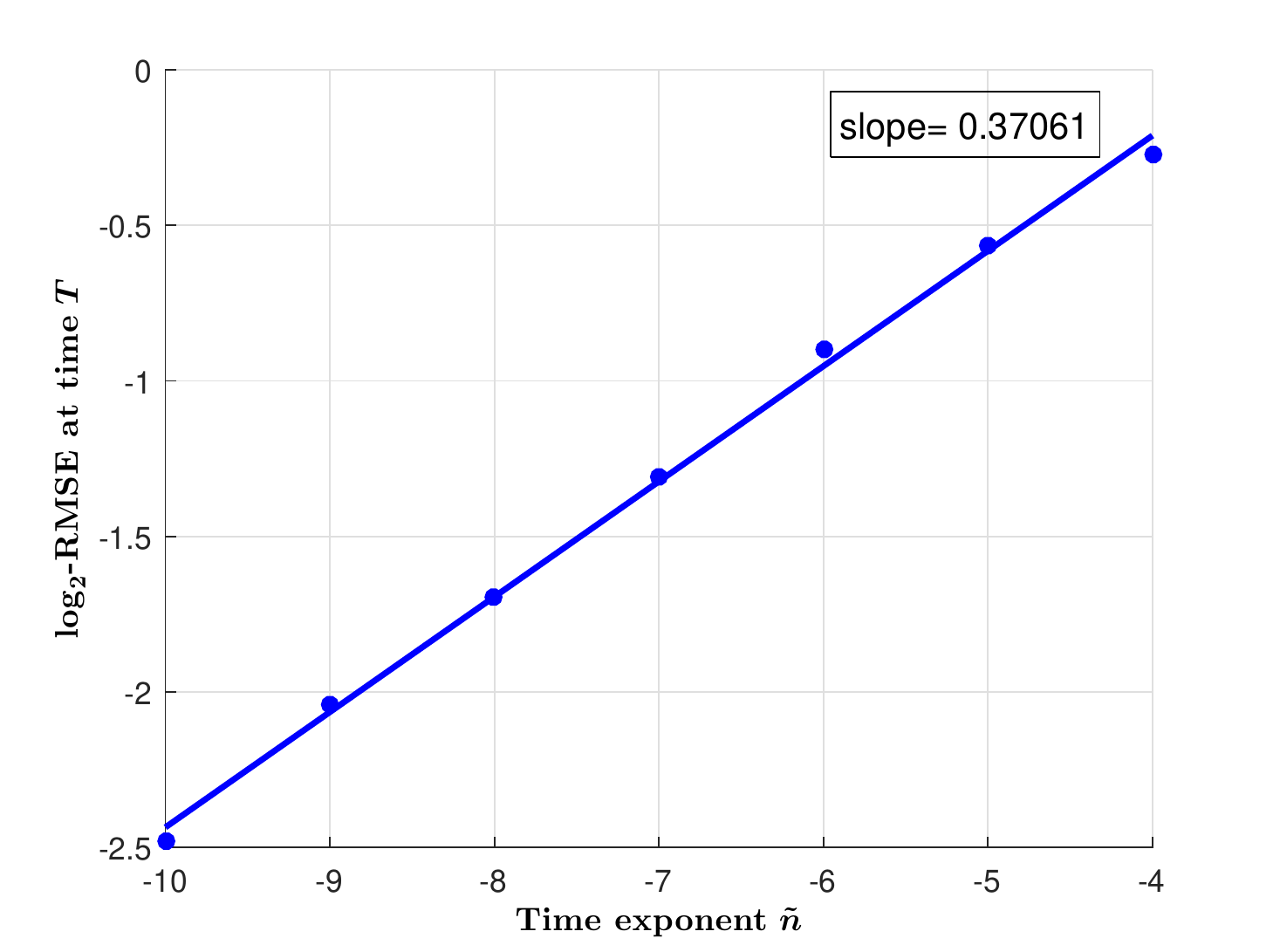}}

	\caption{Euler rates of convergence for elementary\_minus34 for $\xi \in \{0,1\}$}
	
	\label{fig:PlatenEulerRates_elementary_minus34_part1}
\end{figure}

For initial values close to the discontinuity the observed empirical convergence order are in some cases far away from the guarenteed $3/4$, altough the linear regression typically produces stable results, see Figure \ref{fig:PlatenEulerRates_elementary_minus34_part1}(b). A possible explanation for this are again the first drift changes. When starting close to the initial value, the first drift changes seem to be very sensitive to the step-size, which results in rather different trajectories of the Euler scheme.
 
Furthermore, if the initial value is far away from the discontinuity, only very few drift changes occur in the underlying SDE (if at all). Hence, if the step size of the Euler scheme is significantly small, these changes are captured and the error drops drastically. Figure \ref{fig:PlatenEulerRates_elementary_minus34_part1} illustrates this by comparing the regressions for an initial value $\xi=0$ away from the discontinuity in $1.4$ and an initial value $\xi=1$, which is closer to the discontinuity. (The regression has also to deal in Figure \ref{fig:PlatenEulerRates_elementary_minus34_part1}(a) with two different regimes.)
Note that the Euler scheme for \eqref{eq:SDE} is always exact up to the time of the first drift change.

Moreover, for an outward pointing drift coefficient, the Euler scheme and the exact solution coincide with high probability, which explains e.g. the errors close to machine accuracy for the drift coefficient $\textrm{sign}$ and the initial value $\xi=5$.
Note that in this setting, the number of paths with at least one drift change is even zero over all saved $10^4$ solution paths.

To explain this phenomenon, consider again the SDE
\begin{align} dX_t=  \left( \alpha_1  \cdot \ind_{(-\infty,0)} (X_t)  + \alpha_2  \cdot \ind_{[0,\infty)} (X_t)   \right) dt + dW_t, \quad t \geq 0,  \qquad X_0=\xi, \label{test_illustrate_2} \end{align}
with $\alpha_1<0<\alpha_2$. An application of formula (5.13) in chapter 3.5.C in \cite{ks} gives
\begin{align} \label{exit_prob_2}
\mathbb{P} \left( \inf_{t \geq 0} |X_t| >0 \right) &= 1-e^{2 \alpha_1 \xi^{-} -2 \alpha_2 \xi^{+}}, \qquad \xi \neq 0.
\end{align}
Note that an initial value $\xi \neq 0$ is not a restriction as we analyze the case of an initial value far away from the discontinuity. So, for drift values $-\alpha_1=\alpha_2=1$, and an initial value $\xi=5$, the Euler scheme is exact with a probability of at least $1-e^{-10}$ $\approx 0.99995460 \ldots $ 

To summarize: Standard Monte Carlo simulations for testing convergence rates seem to be unreliable in the case of outward pointing coefficients. 
No stable asymptotic regime seems to be reached by our estimators. Smaller stepsizes or a larger Monte-Carlo sample might be a remedy for this problem, similar to \cite{arnulf} where moment explosions of the Euler scheme for SDEs with superlinear cofficients are observed in a numerically asymptotic setting. But this is beyond the scope of the present manuscript.


\section{The Euler scheme for the Atlas model}\label{sec:finance}


In this section, we will use the Euler scheme to simulate the so-called Atlas model, which is a particular first-order market model \cite{Banner.2005}. In such models, the asset dynamics depend
on the size (measured in terms of market capitalization) of the corresponding firm, which results in an SDE model with discontinuous coefficients.

\medskip

\subsection{First-order market models} \label{subsec:FO-model}
A \textit{first-order model} \cite{Banner.2005} is defined as follows:
Let $\gamma, g_1,...,g_d \in \mathbb{R}$ and $\sigma_1,...,\sigma_d \in (0, \infty)$ such that
\begin{align*}
	g_1<0, \quad g_1+g_2<0, \ldots,  \quad g_1+\cdots+g_{d-1}<0, \quad g_1+\cdots+g_d=0.
\end{align*}
Consider now stocks for which the market capitalizations are given by $X_1,\ldots,X_d$, where the index  $ i \in \{1,2, \ldots, d \}$ indicates the name of the firm, and that follow the dynamics
\begin{align} \label{eq:FOmodel}
	d\log X_i(t) = \gamma_i(t) dt + \sigma_i(t) dW_i(t), \quad t \in [0,\infty),\qquad i=1,\ldots,d.
\end{align}
Here, $W_1,...,W_d$ are independent Brownian motions and the {growth rates} $\gamma_i:[0,\infty) \rightarrow \mathbb{R}$ and {volatilities} $\sigma_i:[0,\infty) \rightarrow (0, \infty)$ are given by
\begin{align} \label{eq:FOmodelParameters}
	\gamma_i(t) & = \gamma + \sum_{k=1}^{d} g_k \ind_{\{r_i(t)=k\}}, \qquad \quad \sigma_i(t) = \sum_{k=1}^{d} \sigma_k \ind_{\{r_i(t)=k\}}.
\end{align}
The ranks $r_i(t)$ for the stock $X_i(t)$ at time $t$ arise from the reverse order-statistics:
\begin{align} \label{eq:ranking}
	\max_{1\leq i\leq d} X_i(t) \defgr X_{(1)}(t) \geq X_{(2)}(t) \geq \cdots \geq X_{(d-1)}(t) \geq X_{(d)}(t) \defgl \min_{1 \leq i \leq d} X_i(t).
\end{align}
Ties in the ranking are resolved by giving the firm with a lower  index $i$ the better ranking. So, in such a model the $k$-th largest firm is assigned a growth rate of $\gamma + g_k$ and a volatility of $\sigma_k$ over the whole time horizon.

\smallskip

According to \cite{Banner.2005}, the simplest among the first-order models is the so-called \textit{Atlas model}, which was introduced in \cite[Ex.\ 5.3.3]{Fernholz.2002}.
Within the setting of \eqref{eq:FOmodel} and \eqref{eq:FOmodelParameters}, choosing
\begin{align} \label{eq:AtlasParameters}
	\gamma = g > 0, \quad g_k = -g, \,\, k=1, \ldots ,d-1, \quad g_d = (d-1)g \quad \text{and} \quad \sigma_i(t) = \sigma > 0, \,\, i=1, \ldots, d,
\end{align}
leads to the Atlas model. Here, only the smallest stock in the market -- called the {Atlas stock} -- has a nonzero but positive growth rate (for its log-dynamics).

By setting $Y_i(t) \defgl \log X_i(t)$, $i=1,...,d$, as well as plugging in the Atlas parameters \eqref{eq:AtlasParameters} in our first-order model \eqref{eq:FOmodel} -- \eqref{eq:FOmodelParameters}, we obtain the Atlas model in compact form as
\begin{align} \label{eq:AtlasModel}
	d Y_i(t) &= (d \cdot g)  \ind_{\{r_i(t)=d\}} dt + \sigma dW_i(t), \quad i=1,...,d.
\end{align}
As stated in \cite[Prop.\ 2.3]{Banner.2005}, the solution of \eqref{eq:AtlasModel} satisfies the ergodic relation
	\begin{align} \label{eq:longTermAtlas}
		 \lim\limits_{T \rightarrow \infty} \frac{1}{T} \int_{0}^{T} \ind_{\{r_i(t)=k\}} dt = \frac{1}{d} \quad \textrm{a.s.}, \qquad i,k=1, \ldots, d, 
	\end{align}
i.e., all stocks in the market asymptotically spent at each rank approximately the same amount of time. Similar  ergodic relations also hold  for general first-order market models.


\subsection{Numerical results}


For simulations of the Atlas and general first-order models one has to rely on discretization schemes such as the Euler method.
In this subsection, we test whether the Euler scheme is able to recover the long time behavior \eqref{eq:longTermAtlas}, i.e.,
whether the discrete occupation rates 
$$   \frac{1}{T} \sum_{\ell=1}^{T/\Delta} \ind_{\{\widehat{r}_i(\ell \Delta)=k\}}, \qquad i,k=1, \ldots, d, $$
where $\widehat{r}_i$ is the discretized counterpart of \eqref{eq:ranking} based on the Euler scheme and $T/\Delta \in \mathbb{N}$, converge to the analytical value.

Here, we consider a three-dimensional model with initial log-capitalizations
${Y(0)=[3.4, 4.1, 5.7]}$ and $\widetilde{Y}(0)=[1.2, 3.5, 10.8]$,  $\gamma=0.1$ as market drift and $\sigma=0.09$ as market volatility\footnote{The market parameters are inspired by parameters from A.\ Banner's (INTECH Investement Technologies LLC, Princeton) presentation on ``Equity Market Stability'' given at the WCMF6 conference, Santa Barbara, 2014.}.
Table  \ref{tab:LTrankings} presents the discrete occupation rates (averaged over $M=10^3$ repetitions) for $\Delta=2^{-14}$ and diffe\-rent values of $T$ as well as the sum of the squared deviations from the analytical asymptotic occupation rate. 
As hoped, the discrete occupation rates converge to the analytical asymptotic occupation rate of $1/d=1/3$ with increasing time horizon.

Furthermore, results suggest that less varying initial capitalizations imply that the numerical values are closer to the analytical result already for shorter time horizons, which coincides with the intuitive understanding.
We also simulated the above scenarios with $\Delta=2^{-10}$ instead of $\Delta=2^{-14}$: all occupation times where equal with an accuracy of four digits and one third of the $90$ occupation rates differed in the fifth digit. This suggests that -- as soon as the step size is small enough -- a further refinement of the step size is no longer beneficial and the crucial simulation parameter is $T$, the endpoint of the considered time horizon.

\bigskip

\begin{table}[htbp]
	\centering
	\caption{Discrete occupation rates for the discretized Atlas model}
	\begin{tabular}{|c|c|l|r|r|r|c|c|c|}
		\toprule
		& \textbf{T} &       & \multicolumn{1}{l|}{\textbf{Firm 1}} & \multicolumn{1}{l|}{\textbf{Firm 2}} & \multicolumn{1}{l|}{\textbf{Firm 3}} & \multicolumn{3}{c|}{\textbf{Quadratic deviations}} \\
		\midrule
		\multirow{14}[10]{*}{$Y(0)$} & \multirow{3}[2]{*}{\textbf{100}} & \textbf{Rank 1} & 0.2911 & 0.2895 & 0.4194 & \multirow{3}[2]{*}{0.0030} & \multirow{3}[2]{*}{0.0031} & \multirow{3}[2]{*}{0.0111} \\
		&       & \textbf{Rank 2} & 0.3425 & 0.3662 & 0.2913 &       &       &  \\
		&       & \textbf{Rank 3} & 0.3664 & 0.3443 & 0.2892 &       &       &  \\
		\cmidrule{2-9}
		& \multirow{3}[1]{*}{\textbf{250}} & \textbf{Rank 1} & 0.3156 & 0.3161 & 0.3683 & \multirow{3}[1]{*}{0.0005} & \multirow{3}[1]{*}{0.0005} & \multirow{3}[1]{*}{0.0018} \\
		&       & \textbf{Rank 2} & 0.3375 & 0.3463 & 0.3162 &       &       &  \\
		&       & \textbf{Rank 3} & 0.3469 & 0.3376 & 0.3155 &       &       &  \\
		\cmidrule{2-9}
		& \multirow{3}[1]{*}{\textbf{500}} & \textbf{Rank 1} & 0.3238 & 0.3246 & 0.3517 & \multirow{3}[1]{*}{0.0001} & \multirow{3}[1]{*}{0.0001} & \multirow{3}[1]{*}{0.0005} \\
		&       & \textbf{Rank 2} & 0.3357 & 0.3400 & 0.3243 &       &       &  \\
		&       & \textbf{Rank 3} & 0.3406 & 0.3354 & 0.3240 &       &       &  \\
		\cmidrule{2-9}
		& \multirow{3}[1]{*}{\textbf{750}} & \textbf{Rank 1} & 0.3273 & 0.3273 & 0.3454 & \multirow{3}[1]{*}{0.0001} & \multirow{3}[1]{*}{0.0001} & \multirow{3}[1]{*}{0.0002} \\
		&       & \textbf{Rank 2} & 0.3347 & 0.3379 & 0.3274 &       &       &  \\
		&       & \textbf{Rank 3} & 0.3380 & 0.3348 & 0.3272 &       &       &  \\
		\cmidrule{2-9}
		& \multirow{3}[1]{*}{\textbf{1000}} & \textbf{Rank 1} & 0.3288 & 0.3287 & 0.3425 & \multirow{3}[1]{*}{0.0000} & \multirow{3}[1]{*}{0.0000} & \multirow{3}[1]{*}{0.0001} \\
		&       & \textbf{Rank 2} & 0.3344 & 0.3368 & 0.3288 &       &       &  \\
		&       & \textbf{Rank 3} & 0.3368 & 0.3345 & 0.3287 &       &       &  \\
		\midrule
		\multicolumn{1}{|r}{} & \multicolumn{1}{r}{} & \multicolumn{1}{r}{} & \multicolumn{1}{r}{} & \multicolumn{1}{r}{} & \multicolumn{1}{r}{} & \multicolumn{1}{r}{} & \multicolumn{1}{r}{} &  \\
		\midrule
		& \textbf{T} &       & \multicolumn{1}{l|}{\textbf{Firm 1}} & \multicolumn{1}{l|}{\textbf{Firm 2}} & \multicolumn{1}{l|}{\textbf{Firm 3}} & \multicolumn{3}{c|}{\textbf{Quadratic deviations}} \\
		\midrule
		\multirow{19}[10]{*}{$\widetilde{Y}(0)$} & \multirow{3}[2]{*}{\textbf{100}} & \textbf{Rank 1} & 0.1464 & 0.1447 & 0.7089 & \multirow{3}[2]{*}{0.0554} & \multirow{3}[2]{*}{0.0562} & \multirow{3}[2]{*}{0.2116} \\
		&       & \textbf{Rank 2} & 0.3883 & 0.4654 & 0.1463 &       &       &  \\
		&       & \textbf{Rank 3} & 0.4653 & 0.3899 & 0.1448 &       &       &  \\
		\cmidrule{2-9}
		& \multirow{3}[1]{*}{\textbf{250}} & \textbf{Rank 1} & 0.2581 & 0.2579 & 0.4840 & \multirow{3}[1]{*}{0.0090} & \multirow{3}[1]{*}{0.0090} & \multirow{3}[1]{*}{0.0341} \\
		&       & \textbf{Rank 2} & 0.3555 & 0.3863 & 0.2582 &       &       &  \\
		&       & \textbf{Rank 3} & 0.3864 & 0.3558 & 0.2577 &       &       &  \\
		\cmidrule{2-9}
		& \multirow{3}[1]{*}{\textbf{500}} & \textbf{Rank 1} & 0.2949 & 0.2956 & 0.4096 & \multirow{3}[1]{*}{0.0023} & \multirow{3}[1]{*}{0.0023} & \multirow{3}[1]{*}{0.0087} \\
		&       & \textbf{Rank 2} & 0.3448 & 0.3599 & 0.2953 &       &       &  \\
		&       & \textbf{Rank 3} & 0.3603 & 0.3446 & 0.2951 &       &       &  \\
		\cmidrule{2-9}
		& \multirow{3}[1]{*}{\textbf{750}} & \textbf{Rank 1} & 0.3082 & 0.3079 & 0.3840 & \multirow{3}[1]{*}{0.0010} & \multirow{3}[1]{*}{0.0010} & \multirow{3}[1]{*}{0.0038} \\
		&       & \textbf{Rank 2} & 0.3407 & 0.3513 & 0.3081 &       &       &  \\
		&       & \textbf{Rank 3} & 0.3511 & 0.3409 & 0.3080 &       &       &  \\
		\cmidrule{2-9}
		& \multirow{3}[1]{*}{\textbf{1000}} & \textbf{Rank 1} & 0.3143 & 0.3142 & 0.3715 & \multirow{3}[1]{*}{0.0006} & \multirow{3}[1]{*}{0.0006} & \multirow{3}[1]{*}{0.0022} \\
		&       & \textbf{Rank 2} & 0.3391 & 0.3467 & 0.3143 &       &       &  \\
		&       & \textbf{Rank 3} & 0.3467 & 0.3391 & 0.3143 &       &       &  \\
		\bottomrule
	\end{tabular}%
	\label{tab:LTrankings}%
\end{table}%

\bigskip

\section{Conclusion and Outlook}
We have seen that the numerical approximation of solutions of SDEs with discontinuous drift coefficients is a challenging task, where several particularities arise. We were able to identify two main classes of discontinuous drift coefficients: outward and inward pointing drift coefficients. For the latter class, we analyzed stability properties.
It turned out that the main difficulty in measuring the empirical convergence rates is how to appropriately capture drift changes. For inward pointing coefficients, we obtained stable estimates, which are in accordance with the theoretical results. For outward pointing cases, the estimates seem to be unreliable, no stabilizing asymptotic regime seems to be reached for the estimates.
 We tested two higher-order numerical schemes, that are frequently used in a setting where coefficients are sufficiently smooth. However, both schemes did not lead to an improved  behavior.

\bigskip



\section*{Acknowledgment}
This work was supported by the DFG grant No.\ GO 1920/4-1. 
Part of this work was carried out while A. Neuenkirch was visiting 
the Facultad de Matem\'aticas de la Universidad de Sevilla; A. Neuenkirch whishes to thank the
Dpto. Ecuaciones Diferenciales y An\'alisis Num\'erico 
for its hospitality and support.

The publication of this article was funded by the Ministry of Science, Research and the Arts Baden-W\"urttemberg and the University of Mannheim.


\bibliography{literature_rankBased_submit}
\bibliographystyle{siam}	


\end{document}